\documentclass[11pt]{amsart}

\usepackage[square,compress,comma, numbers,sort]{natbib}
\usepackage[colorlinks=true, citecolor=red, linkcolor=blue]{hyperref}
\usepackage{amsfonts,mathtools}

\allowdisplaybreaks[4]

\usepackage{amsmath}
\usepackage{bbm}
\usepackage{amssymb}
\usepackage{mathtools}
\usepackage[dvipsnames]{xcolor}

\usepackage{mathrsfs}
\usepackage[scr=rsfs,cal=boondox]{mathalfa}

\usepackage{color}

\definecolor{c20}{rgb}{0.,0.7,0.}
\definecolor{c30}{rgb}{0.,0.,1.}
\definecolor{c40}{rgb}{1,0.1,0.7}
\definecolor{c50}{rgb}{1,0,0}
\definecolor{c60}{rgb}{1,0.9,0.1}

\newcommand{\abs}[1]{\left\lvert #1 \right\rvert}

\newcommand{\sprod}[1]{\langle#1\rangle}

\newcommand{\E}[1]{\mathbb{E}\left\{ #1\right\}}

\newcommand{\pk}[1]{\mathbb{P} \left\{ #1 \right \} }

\newcommand{\R}{\mathbb{R}}

\newcommand{\N}{\mathbb{N}}

\newcommand{\BQN}{\begin{eqnarray}}
\newcommand{\EQN}{\end{eqnarray}}
\newcommand{\BQNY}{\begin{eqnarray*}}
\newcommand{\EQNY}{\end{eqnarray*}}

\newcommand{\BS}{\begin{sat}}
\newcommand{\ES}{\end{sat}}
\newcommand{\BT}{\begin{theo}}
\newcommand{\ET}{\end{theo}}
\newcommand{\BK}{\begin{korr}}
\newcommand{\EK}{\end{korr}}

\newcommand{\BD}{\begin{de}}
\newcommand{\ED}{\end{de}}

\newcommand{\BIT}{\begin{itemize}}
\newcommand{\EIT}{\end{itemize}}
\newcommand{\BDI}{\begin{description}}
\newcommand{\EDI}{\end{description}}

\newcommand{\BRM}{\begin{remarks}}
\newcommand{\ERM}{\end{remarks}}

\newcommand{\BEL}{\begin{lem}}
\newcommand{\EEL}{\end{lem}}

\newtheorem{theo}{Theorem}[section]
\newtheorem{sat}[theo]{Proposition}
\newtheorem{de}[theo]{Definition}
\newtheorem{lem}[theo]{Lemma}

\newtheorem{example}[theo]{Example}
\newtheorem{korr}[theo]{Corollary}
\newtheorem{remark}[theo]{Remark}
\newtheorem{remarks}[theo]{Remarks}

\newcommand{\nelem}[1]{{Lemma \ref{#1}}}

\newcommand{\netheo}[1]{{Theorem \ref{#1}}}

\newcommand{\prooftheo}[1]{ \textsc{\bf Proof of Theorem} \ref{#1}:}

\newcommand{\prooflem}[1]{\textsc{\bf Proof of Lemma} \ref{#1}:}
\newcommand{\proofrem}[1]{\textsc{\bf Proof of Remark} \ref{#1}:}

\newcommand{\COM}[1]{}

\def\td{\text{\rm d}}

\newcommand{\QED}{\hfill $\Box$}

\newcommand{\kb}[1]{\boldsymbol{#1}}
\newcommand{\vk}[1]{\kb{#1}}

    \title[Functional limit theorems for Gaussian-fed queueing networks]{Functional limit theorems for Gaussian-fed \\ queueing network in light and heavy traffic}

\author{Nikolai Kriukov}
\address{Nikolai Kriukov, Korteweg-de Vries Institute\\
	University of Amsterdam\\
	P.O. Box  94248, 1090 GE  Amsterdam, The Netherlands
}
\email{n.kriukov@uva.nl}

\author{Krzysztof D\c{e}bicki}
\address{Krzysztof D\c{e}bicki, Mathematical Institute\\
        Wroc\l{}aw University\\
        pl.\ Grunwaldzki 2/4, 50-384 Wroc\l{}aw, Poland}
\email{krzysztof.debicki@math.uni.wroc.pl}

\author{Michel Mandjes}
\address{Michel Mandjes, Mathematical Institute\\
        Leiden University\\
        P.O. Box 9512, 2300 RA Leiden, The Netherlands}
\email{m.r.h.mandjes@math.leidenuniv.nl}

\begin{document}

\maketitle

\begin{abstract}
   {We} consider a queueing network operating under a strictly upper-triangular routing matrix with per
   column at most one non-negative entry. The root node is fed by a Gaussian process with stationary increments.
   Our aim is to characterize the distribution of the multivariate stationary workload process under a specific scaling
   of the queue's service rates. In the main results of this paper we identify,  under mild conditions on the standard
   deviation function of the driving Gaussian process, in both light and heavy traffic {parameterization}, the limiting
   law of an appropriately scaled version (in both time and space) of the joint stationary workload process.
   In particular, we develop conditions under which specific {queueing processes} of the network effectively decouple, i.e.,
   become independent in the limiting regime.

    \vspace{3mm}


\vspace{3mm}

\noindent
{\sc Affiliations.} NK is affiliated with the Korteweg-de Vries Institute for Mathematics, University of Amsterdam, Science Park 904, 1098 XH Amsterdam, The Netherlands.

\noindent
KD is affiliated with Mathematical Institute, Wroc\l{}aw University, pl.\ Grunwaldzki 2/4, 50-384 Wroc\l{}aw, Poland.

\noindent
MM is affiliated with with the Mathematical Institute, Leiden University, P.O. Box 9512,
2300 RA Leiden,
The Netherlands. MM is also affiliated with Korteweg-de Vries Institute for Mathematics, University of Amsterdam, Amsterdam, The Netherlands; E{\sc urandom}, Eindhoven University of Technology, Eindhoven, The Netherlands; Amsterdam Business School, Faculty of Economics and Business, University of Amsterdam, Amsterdam, The Netherlands.

\vspace{3mm}

\noindent
{\sc Acknowledgments.}
NK and MM are supported by the European Union’s Horizon 2020 research and innovation programme under the Marie Sk\l{}odowska-Curie grant agreement no.\ 945045, and by the NWO Gravitation project {\tiny NETWORKS} under grant agreement no.\ 024.002.003. \includegraphics[height=1em]{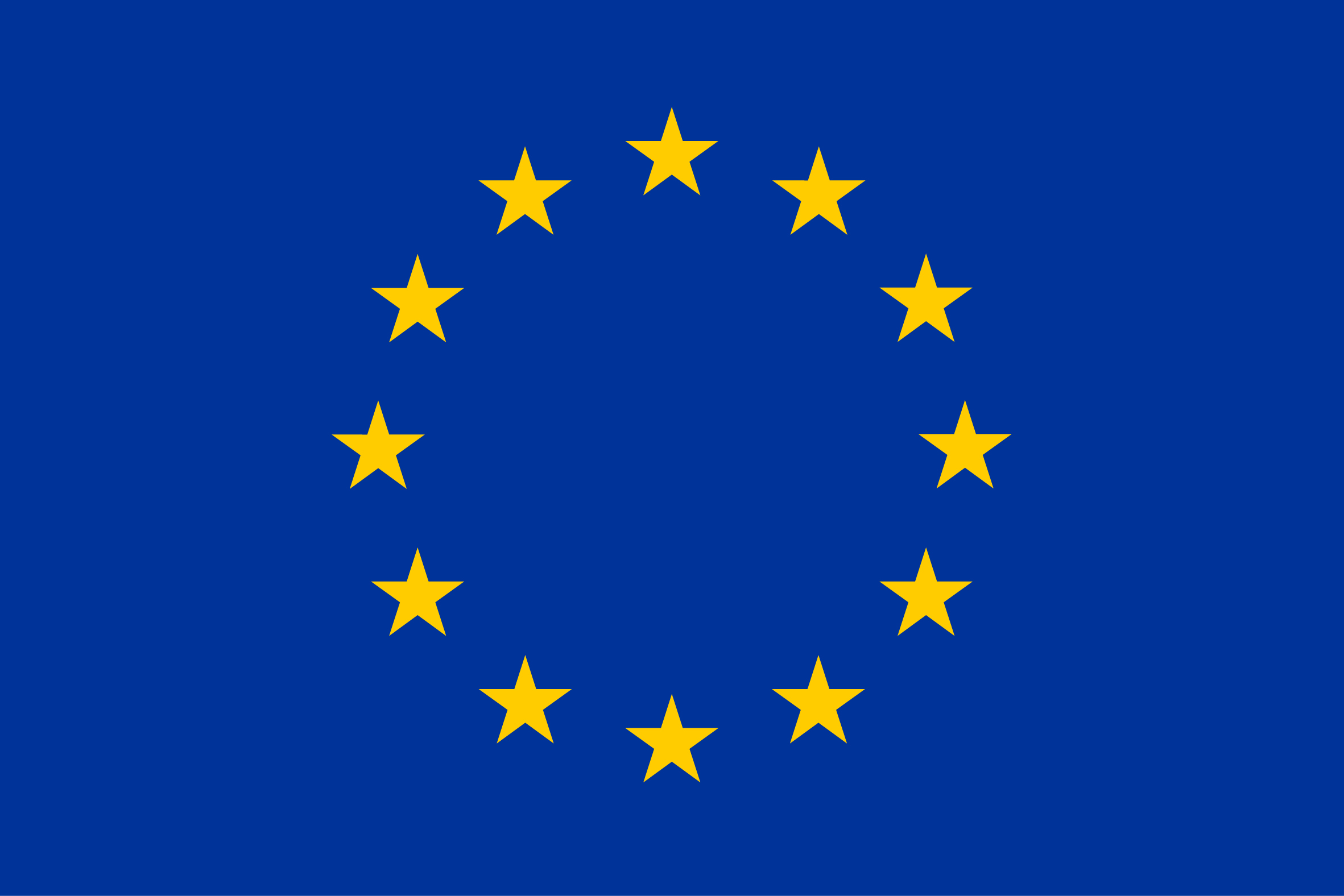}

\end{abstract}

\newpage

\section{Introduction}\label{Section:Introduction}

Consider the following classical model that describes the evolution of a {\it storage system}, also sometimes referred to as the {\it workload} of a {\it queue}. There is a buffer, of unlimited capacity, which is fed by an external input process $J(\cdot)$ that we throughout assume to have stationary increments, and which is emptied at a constant output rate $r>0$ (to be interpreted as the {\it maximal} rate at which work exits the buffer). The storage level, or workload, at time $t$ is given in terms of the model primitives $J(\cdot)$ and $r$ via the representation
\begin{align*}
    Q(t) = Q(0) + J(t) - rt +\max\biggl(0, \sup_{0<s<t}\bigl(-Q(0) - J(s)+rs\bigr)\biggr).
\end{align*}
Importantly, this formalism does not require that the process $J(\cdot)$ be increasing.
In the sequel, we focus on the system being in stationarity at time $0$,
where we impose the stability condition $r>{\mathbb E}[J(1)]$.
In this case we can write the workload process as the reflected version of the `free process' $(J(t)-rt)_{t\in{\mathbb R}}$:
\begin{align*}
    Q(t) = \sup_{-\infty <s < t}\bigl(J(t) - J(s) - r(t-s)\bigr);
\end{align*}
see e.g.\ \cite[\S 2.4]{DeM15}.

While single queues have been extensively studied, significantly less is known about their network counterparts. The simplest instance of such a network is the two-node {\it tandem} queue, which consists of two concatenated queues where the output of the first queue serves as the input for the second. For an $n$-node tandem queue, denoting by $\vk r = (r_1,\ldots,r_n)^\top$ the respective output rates, there is an elementary way to describe the multivariate stationary workload process $\vk Q(\cdot) = (Q_1(\cdot),\ldots, Q_n(\cdot))^{\top}$; see e.g.\ \cite[Remark 3.17]{Mandjes_2005}, and in addition,
for more background on multivariate storage processes, \cite[Chs.\ XII-XIII]{DeM15} or \cite{Rob13}. Here we assume $r_1>\ldots>r_n$, as if for a pair of subsequent queues $i$ and $i+1$ this ordering would be violated, the $(i+1)$-st queue would remain empty and hence can be left out. Based on the fact that \[
\sum_{j=1}^i Q_i(t) = \sup_{-\infty < s < t}\bigl(J(t)-J(s) - r_i(t-s)\bigr),\]
we have, for any $t\geqslant 0$ and any $i\in\{2,\ldots,n\}$, under the assumption  $\min_i r_i> {\mathbb E}[J(1)]$,
\begin{align*}
    Q_i(t)  =  \sup_{-\infty < s < t}\bigl(J(t)-J(s) - r_i(t-s)\bigr) -  \sup_{-\infty < s < t}\bigl(J(t)-J(s)  - r_{i-1}(t-s)\bigr).
\end{align*}

A next step is to extend the above class of tandem queues to more general $n$-node {\it queueing networks}. A convenient class of networks for this purpose includes those in which the output of any queue can be \textit{split} according to a fixed fraction, with the resulting individual output processes either feeding into subsequent queues or leaving the network, and in which there are \textit{no loops}; see the formalism studied in, for example, \cite{DDR}. Crucially, this class of networks does not allow streams to \textit{merge}.
With $p_{ij}$ denoting the fraction of the $i$-th node which serves as the input of $j$-th node,  the routing is encoded via the {\it routing matrix} $P = \{p_{ij}\}_{i,j=1}^{n}$. As follows from the assumptions imposed, for the networks that we consider $P$ is strictly upper-triangular and in addition each of its columns contains at most one positive element. In the specific case of the $n$-node tandem queue we have $p_{i,i+1}=1$  for $i\in\{1,\ldots,n-1\}$ and $p_{ij}=0$ otherwise, but note that the class of networks we introduced also encompasses tandem queues in which part of a queue's output instantly leaves  the network (i.e., $p_{i,i+1}\in(0,1)$  for $i\in\{1,\ldots,n-1\}$ while still all other $p_{ij}=0$).

The main advantage of the class of queueing networks that we introduced in the previous paragraph, is that they can be represented in a compact format, as follows.
Denote, as before, by $\vk r = (r_1,\ldots,r_n)^{\top}$ the output rates of the $n$ nodes, and by $J(\cdot)$ the external input process feeding into the first node. Then  the  stationary workload process $\vk Q(\cdot) = (Q_1(\cdot),\ldots,Q_n(\cdot))^{\top}$ is given through \cite{DDR,DeM15}
\begin{align*}
    \vk Q(\cdot) = (I - P^{\top})\overline{\vk X}(\cdot),
\end{align*}
where the process $\overline{\vk X}(\cdot) = \bigl(\overline{X}_{1}(\cdot),\ldots,\overline{X}_{n}(\cdot)\bigr)^{\top}$ is defined for any $t\in[0,\infty)$ as
\begin{align}
    \overline{X}_{i}(t) := \sup_{-\infty < s < t}\bigl(\mathcal{C}_i\bigl(J(t)-J(s)\bigr) - r_i(t-s)\bigr),\:\quad
\mathcal{C}_i := (I-P^{\top})^{-1}_{i,1}.\label{eq:Ci}
\end{align}
Conditions (in terms of $P$ and $\vk r$) must be imposed to ensure the existence of a stationary workload process. Sufficient conditions, that apply in our specific setting, are provided in Section \ref{Section:Model}, but see \cite{DDR} and \cite[Ch.\ XIII]{DeM15} for more background.

Remarkably, queueing networks that fall outside the class described above — such as those with merging streams or loops — are inherently more difficult to analyze. Specifically, in these cases one cannot express the components of $\vk Q(t)$ as a linear combination of the entries of $\overline{\vk X}(t)$. We do not address these more general networks in this paper.

\smallskip

\noindent {\it Scope \& contribution ---}  In this paper we focus on queueing networks of the type described above,
where the driving process $J(\cdot)$ is assumed to be a {\it Gaussian process} with stationary increments.
The use of Gaussian input is motivated by a central-limit type argument:
large aggregates of weakly dependent contributions tend to exhibit Gaussian behavior
\cite{TWS97, Ste02, mandjes2007, Meent2006GaussianTE}.
{From the point of view of stochastic modelling,} the class of Gaussian processes is particularly
appealing because it encompasses a broad range of correlation structures.
Specifically, it includes {\it fractional Brownian motions,  self-similar processes that are},
choosing the Hurst parameter larger than $\frac{1}{2}$, long-range dependent.
{Supported by extensive measurement studies \cite{Cro95, Lel94}, this class of Gaussian processes has become central to traffic modeling in modern communication networks, as argued in works such as \cite{Nor94, mandjes2007}.}

In the single-queue context, the stationary storage level has been studied in detail; see e.g.\ the overview on
large-buffer asymptotics in \cite[\S 5.6]{mandjes2007}. Such a probabilistic analysis of multi-node networks is still lacking,
however, mainly due to the complex nature of the shaping effects that traffic streams undergo as they pass through the network
nodes; for partial results on tandem queues we refer to e.g.\ \cite{Mandjes_2005,manner}.
It is the primary objective of this paper to shed light on this,
where two specific asymptotic regimes are considered, viz.:
\begin{itemize}
\item[$\circ$] {\it light traffic}: for all nodes $i\in\{1,\ldots,n\}$ the output rates $r_i$ are such that $r_i=r_i(u)\to \infty$,
\item [$\circ$] {\it heavy traffic}: for all nodes $i\in\{1,\ldots,n\}$ the output rates $r_i$ are such that $r_i=r_i(u)\to 0$,
\end{itemize}
in the regime where $u\to\infty$, where $u$ is a scaling parameter.
Both light- and heavy-traffic approximations are essential in the asymptotic analysis of queueing systems,
as they provide insightful expressions
for performance measures that are typically highly complex under non-asymptotic conditions.
We refer to e.g.\  \cite{Kin, KBZ, BC, RS, whitt} for related works on heavy-traffic
approximations 
of various non-Gaussian queueing models, and to e.g.\ \cite{DaR, BRS, Sig}
for light-traffic approximations in the setting of conventional queueing models.

{For the Gaussian-fed single queue analyzed in \cite{Gaussianqueues},} it was shown that by imposing a regularly varying
condition on the standard deviation function of the Gaussian input process $J(\cdot)$,
the stationary workload process of a single queue,
after an appropriate rescaling, converges in distribution to a storage system driven by fractional Brownian motion,
in both light and heavy traffic regimes.
The primary goal of this paper is to extend this result to the significantly broader class of queueing networks defined above.

Our {work} contains two main results. \netheo{main_light_traffic} focuses on the light-traffic situation,
establishing a functional limit theorem for the multivariate stationary workload processes ${\boldsymbol Q}(\cdot)$,
under mild conditions on the standard deviation function pertaining to the input process $J(\cdot)$.
\netheo{main_heavy_traffic} is the counterpart of \netheo{main_light_traffic}, but then for the heavy-traffic scenario. 
Theorems~\ref{main_light_traffic} and~\ref{main_heavy_traffic} provide explicit conditions under which individual queues effectively `decouple' — that is, become independent as $u \to \infty$. In the heavy-traffic regime, these conditions align with those recently found in \cite{DAI} for generalized Jackson networks, leading to an asymptotic product-form structure. A key difference between the results in \cite{DAI} and ours is that \cite[Thm. 3.1]{DAI} focuses solely on stationary workloads, whereas we aim to derive limit theorems at the process level. Additionally, while \cite{DAI} imposes no assumptions on the routing matrix, our framework allows for the possibility of workloads at certain stations remaining dependent. In contrast, the setup in \cite{DAI} considers only a regime where {\it all} queues decouple.

{At the methodological level, we apply the theory of {\it stochastic-process limits} \cite{billingsley2013convergence, whitt} to analyze heavy- and light-traffic regimes, which facilitates the establishment of {\it functional limit theorems}. Over the past few decades, various tools have been developed for proving weak convergence at the process level  in the context of Gaussian processes, with continuous mapping theorems playing a key role. However, for the purposes of this paper, these techniques required careful adaptation and integration. In particular, we highlight the challenges arising from the fact that the queueing processes considered here are defined by applying the supremum functional over the left half-line, i.e., $(-\infty, t]$ (see (\ref{eq:Ci})). Due to the discontinuity of the supremum functional in $C((-\infty, T] \to {\mathbb R}^n)$, a precise analysis is necessary to justify passing to the limit.}

\smallskip

\noindent {\it Notation ---} Throughout the paper we consistently use the following {notation}.
For any {natural number} $d$, we write $\sprod{d}$ for the set $\{1,\ldots,d\}$. All vectors are written in bold and, vice versa, all bold symbols refer to vectors. All entries of a vector share the same letter, but do not use bold font and have a subscript which represents the index of the entry. For example if $\vk b\in\R^d$, then $\vk b = (b_1,\ldots,b_d)^{\top}$. If a vector has a subscript itself, for each element of this vector we first write a subscript of the vector, and after the subscript of the element. For example, if $\vk b_i\in\R^d$, then $\vk b_i = (b_{i,1},\ldots,b_{i,d})^{\top}$.  We denote by $\vk b\circ \vk c\in\R^d$ the componentwise product of $\vk b\in\R^d$ and $\vk c\in\R^d$, i.e., a vector whose $i$-th component is $b_ic_i.$ We also use the notation $\vk{1} =(1,\ldots,1)^\top\in\R^n$, and $\vk\infty$ and $\vk 0$ have similar meanings.

\smallskip

\noindent{\it Outline ---} This paper is organized as follows. Section \ref{Section:Model} formally introduces the model and the assumptions imposed. In Section \ref{Section:Main_results} we present the main results of this paper. We state the functional limit theorems for the joint stationary workload process, and present the main lemmas which are needed to prove them. Sections \ref{Section:Proof_of_gauss_conv} and \ref{Section:Proof_of_sup_conv} are dedicated to the proofs of the lemmas presented in Section~\ref{Section:Main_results}. Appendix~\ref{Section:Appendix} contains (relatively technical) proofs of a series of auxiliary results.

\section{Model}\label{Section:Model}
In this section we detail our model, and discuss the assumptions imposed.
A key role is played by the univariate process  $J(\cdot)$, which is a centered Gaussian process with stationary increments,
uniquely characterized through its variance function $\sigma^2(\cdot)$, i.e., $\sigma^2(t) = {\rm Var}\,J(t)$ for $t\geqslant 0$; in the sequel we often work with the corresponding standard deviation function $\sigma(\cdot)$.

The network dynamics {is} defined via the strictly upper-triangular routing matrix
$P= \{p_{ij}\}_{i,j=1}^{n}$, that is such that its $j$-th column contains exactly one positive element for
$j\in\{2,\ldots,n\}$. A deterministic fraction $p_{ij}$ of the output of queue $i$ is fed into queue $j$.
The first node of this network (the `root node') is fed by $J(\cdot)$.
The output rate of the $i$-th node is denoted by functions $r_i = r_i(u)$ that we assume to be monotone in $u$,
where $u$ is a scaling parameter that we let grow large.

We assume that the matrix $P$ and the vector $\vk r(u)$ satisfy the following two properties:
\begin{itemize}
    \item[\textbf{N1}:]{For any $i,j\in\sprod{n}$, if $p_{ij}>0$ then $p_{ij}\,r_{i}(u)>r_j(u)$ for any $u>0$,}
    \item[\textbf{N2}:]{For any $i\in\sprod{n}$ and any $u>0$, {${\mathbb E}[J(1)]\,(I-P^{\top})^{-1}_{i,1}=0<r_i(u)$}.}
\end{itemize}
Regarding ${\bf N2}$, observe that ${\mathbb E}[J(1)]\,(I-P^{\top})^{-1}_{i,1}=0$ since we assumed that the process $J(\cdot)$ is centered.
{Under these assumptions the existence of the stationary workload process is ensured.
Indeed, as pointed out in e.g.\ \cite{DDR} and \cite[Ch.\ XIII]{DeM15}, under these conditions
{the stationary $n$-dimensional queueing process $\vk Q_{u}(\cdot)$, for any given $u>0$, obeys the representation}
\begin{align}
    \vk Q_{u}(\cdot) = (I - P^{\top})\overline{\vk X}_u(\cdot),\label{Q_vk_def}
\end{align}
where $\overline{\vk X}_u(\cdot) = \bigl(\overline{X}_{u,1}(\cdot),\ldots,\overline{X}_{u,n}(\cdot)\bigr)^{\top}$ is
defined by
\begin{align}
    \overline{X}_{u,i}(t) := \sup_{-\infty < s < t}\bigl(\mathcal{C}_i\bigl(J(t)-J(s)\bigr) - r_i(u)(t-s)\bigr),\label{X_bar_def}
\end{align}
for any $t\in{\mathbb R}$ and with ${\mathscr C}_i$ as given in \eqref{eq:Ci}.}

Throughout this paper we consider two limiting regimes: the {\it light-traffic} scenario in which
$\vk r(u)\to \vk \infty$ as $u\to\infty$, and the {\it heavy-traffic} scenario in which $\vk r(u)\to 0$ as $u\to\infty$.
In each of these two regimes we impose specific assumptions on the shape of the standard deviation function underlying
the process $J(\cdot)$. We denote by ${\mathscr R \mathscr V}_i(\lambda)$ the class of regularly varying functions
at $i\in\{0,\infty\}$ with index~$\lambda$; see e.g.\ \cite{bingham}.

\begin{itemize}
    \item[$\circ$]
For the light-traffic case we assume that the standard deviation $\sigma(\cdot)$ of $J(\cdot)$ satisfies:
\begin{itemize}
	\item[\textbf{L1}:]{$\sigma(\cdot)\in\operatorname{\mathscr{RV}}_0(\lambda)$ for some $\lambda\in(0,1)$,}
        \item[\textbf{L2}:]{For some $t_0>0$ there exist constants $\beta\in(0,1)$, $C>0$ such that, for all $t>t_0$,
        \begin{align*}
            \sigma(t)\leqslant Ct^{\beta}.
        \end{align*}}
\end{itemize}
\item[$\circ$] For the heavy-traffic we assume  that the standard deviation $\sigma(\cdot)$ of $J(\cdot)$ satisfies:
\begin{itemize}
	\item[\textbf{H1}:]{$\sigma(\cdot)\in\operatorname{\mathscr{RV}}_{\infty}(\alpha)$ for some $\alpha\in(0,1)$,}
        \item[\textbf{H2}:]{For some $t_0>0$ there exist constants $\beta\in(0,1]$, $C>0$ such that, for all $t<t_0$,
        \begin{align*}
            \sigma(t)\leqslant Ct^{\beta}.
        \end{align*}}
\end{itemize}
\end{itemize}
Observe that these conditions are mild, in that they are fulfilled by a wide range of relevant Gaussian processes with stationary increments; see, e.g.,
the processes discussed in \cite{AMN}.
{For the both regimes, the above conditions imply the existence of a modification of
$J(\cdot)$ which has continuous sample paths a.s. In the rest of this contribution, we consider
only such modifications.}

\begin{remark}\label{rem:H2_L2}\em
Note that the light-traffic and heavy-traffic assumptions are nearly symmetrical, yet they differ subtly.
   Allowing $\beta = 1$ in \textbf{H2} ensures that the heavy-traffic result also applies to the case where the driving process
   $J(\cdot)$ is an integrated Gaussian process, i.e., $J(t) = \int_0^t \eta(s) , \mathrm{d}s$, with $\eta(\cdot)$ a stationary,
   a.s. continuous,   centered Gaussian process (which has a standard deviation function that grows linearly at 0).
   This is not the case for the light-traffic scenario, as \textbf{L1} excludes such integrated processes.
   {We refer to Remark \ref{rem_deg}
   for a more in-depth explanation of this subtle difference; see also \cite{Gaussianqueues}.}
\hfill    ${\Diamond}$
\end{remark}

Bearing in mind the assumptions \textbf{L1} and \textbf{H1}, it is convenient to define, for any $x\in(0,\infty)$ and $t\in\R$, the function
\begin{align}
    \mathcal{F}_x(t) := \frac{\sigma(tx)}{\sigma(x)}.\label{F_cal_def}
\end{align}
Then, Assumption \textbf{L1} entails that, for any fixed $t\in\R$,
\begin{align}
    \lim_{x\downarrow 0}\mathcal{F}_{x}(t) = \abs{t}^{\lambda},\label{RV_statement}
\end{align}
while Assumption \textbf{H1} yields, for any fixed $t\in\R$,
\begin{align}
    \lim_{x\to \infty}\mathcal{F}_{x}(t) = \abs{t}^{\alpha}.\label{F_RV_infty}
\end{align}

As it turns out, the joint convergence results we are targeting require us to impose the following assumption on the relative asymptotic behavior of the functions $r_i(\cdot)$.
\begin{itemize}
\item[\textbf{N3}:]{For any $i,j\in\langle n\rangle$ such that $i\leqslant j$, the following limit exists:
    \begin{align*}
        \mathfrak{r}_{ij}:=\lim_{u\to\infty} \frac{r_j(u)}{r_i(u)}\in[0,1].
    \end{align*}
    }
\end{itemize}
\begin{remark}\em
    The assumption $\mathfrak{r}_{ij}\in[0,1]$ can be regarded as non-restrictive for the following reason. {If a network satisfies assumptions \textbf{N1}, \textbf{N2}, and the limits $\mathfrak{r}_{ij}$ exist, but do not necessarily belong to $[0,1]$, we can always relabel the nodes of this network in such a way that all three assumptions \textbf{N1}, \textbf{N2} and \textbf{N3} are satisfied. We can do this as follows. Consider a queueing network with $n$ nodes, strictly upper-triangular routing matrix $P$ and output rates $r_1(u),\ldots, r_n(u)$ such that it satisfies \textbf{N1}, \textbf{N2}, and the limits $\mathfrak{r}_{ij}$ presented in \textbf{N3} exist for all $i,j\in\sprod{n}$, $i\leqslant j$.} {For such a network we can always find} a permutation ${\boldsymbol \pi}:\sprod{n}\to\sprod{n}$ such that $i<j$ implies $\mathfrak{r}_{\pi(i),\pi(j)}\in[0,1]$, and $\mathfrak{r}_{\pi(i),\pi(j)} = 1$ implies $\pi(i)<\pi(j)$ if $i<j$. Now, instead of our actual queuing network we  consider an alternative one: it is fed by the same external input process $J(\cdot)$, but the routing matrix and output rates become $P^{\prime}=\{p'_{ij}\}_{i,j=1}^n$ and $\vk{r}^{\prime}(u)$, respectively, which are defined via \[p_{ij}':= p_{\pi(i),\pi(j)},\quad r^{\prime}_i(u) :=r_{\pi(i)}(u).\]

    We proceed by checking that this defines a queuing network that is contained in the class we introduced above. This in the first place means that we have to check whether $P^{\prime}$ is strictly upper-triangular. Assume that we have $i,j\in\sprod{n}$ such that $i<j$ and $p^{\prime}_{ji}>0$. Then, on one hand, we know that $p_{\pi(j),\pi(i)}>0$, hence due to \textbf{N1} for any $u>0$
    \begin{align*}
        r_{\pi(j)}(u)>r_{\pi(i)}(u),
    \end{align*}
    meaning that $\mathfrak{r}_{\pi(i),\pi(j)}\geqslant 1$. On the other hand, as $i<j$, we have $\mathfrak{r}_{\pi(i),\pi(j)}\geqslant 1$.
    Hence, $\mathfrak{r}_{\pi(i),\pi(j)}= 1$ implies that $\pi(i)<\pi(j)$, which contradicts the assumption that $p_{\pi(j),\pi(i)}>0$, as $P$ is strictly upper-triangular.

    {Hence, $P^{\prime}$ is upper-triangular,} meaning that the constructed queueing network is in the class of networks considered. {From the construction of the new network we can see that it satisfies \textbf{N1}, \textbf{N2} and \textbf{N3}, and, denoting by $\vk Q^{\prime}_{u}(\cdot)$ its stationary workload process, the following equation holds in distribution in $C([0,\infty)\to\R^n)$:  for any $u>0$,}
    \begin{align*}
        \left(Q^{\prime}_{u,1}(\cdot),\ldots Q^{\prime}_{u,n}(\cdot)\right)^{\top} = \left(Q_{u,\pi(1)}(\cdot),\ldots Q_{u,\pi(n)}(\cdot)\right)^{\top}.
    \end{align*}
    Consequently, we can consider the newly constructed network rather than the actual one.\COM{, bearing in mind that the stationary workload process of the constructed network is (up to the permutation ${\boldsymbol \pi}$)  equivalent to the stationary workload process of the actual network, i.e.}
    \hfill    ${\Diamond}$
\end{remark}

A result presented in \cite{Gaussianqueues} provides that each process $\overline{X}_{u,i}(\cdot)$, after an appropriate rescaling of time and space, has a distributional limit in $C([0,\infty)\to\R)$. To define  the correct  scaling, we introduce the function $\delta(\cdot)$
as the solution of the following equation (where it exists):
\begin{align}
    \frac{x\delta(x)}{\sigma(\delta(x))} = 1. \label{delta_def}
\end{align}
Then define {the processes $\overline{\vk X}^{\delta}_u(\cdot)$ and $\vk Q^{\delta}_u(\cdot)$, for any $t\in[0,\infty)$, by
\begin{align}
    \overline{\vk X}_u^{\delta}(t)&:=\left(\frac{\overline{X}_{u,1}(\delta(r_1(u))t)}{\sigma(\delta(r_1(u)))},\ldots,\frac{\overline{X}_{u,n}(\delta(r_n(u))t)}{\sigma(\delta(r_n(u)))}\right)^{\top},\label{X_delta_def}\\
    \vk Q^{\delta}_u(t) &:= (I - P^{\top})\overline{\vk X}_u^{\delta}(t) = \left( \frac{Q_{u,1}\bigl(\delta(r_1(u))t\bigr)}{\sigma\bigl(\delta(r_1(u))\bigr)},\cdots, \frac{Q_{u,n}\bigl(\delta(r_n(u))t\bigr)}{\sigma\bigl(\delta(r_n(u))\bigr)} \right)^{\top}.\notag
\end{align}}
\begin{remark}\em In general, the function $\delta(\cdot)$ may not be defined for some $x\in\R$. However, as argued in \cite{Gaussianqueues}, assumption \textbf{L1} implies that it is well-defined in some vicinity of zero, and
\begin{align}
    \delta(\cdot) \in \mathcal{RV}_{\infty}\bigl(1/(\lambda - 1)\bigr).\label{delta_rv}
\end{align}
Similarly, assumption \textbf{H1} entails that the function $\delta(\cdot)$ is well-defined at infinity, and
\begin{align*}
    \delta(\cdot) \in \mathcal{RV}_{0}\bigl(1/(\alpha - 1)\bigr).
\end{align*}
These properties will be extensively used throughout our derivations.
\hfill    ${\Diamond}$
\end{remark}

The results derived in \cite{Gaussianqueues} show that, under the appropriate assumptions, each component $\overline{X}^{\delta}_{u,i}(\cdot)$ has a limit in $C([0, \infty) \to \mathbb{R})$ as $u \to \infty$ in both the light- and heavy-traffic scenario.
{In Section \ref{Section:Main_results}, by appealing to \eqref{Q_vk_def},
we derive} a functional limit for the full multivariate process $\vk Q^{\delta}_u(\cdot)$ in $C([0, \infty) \to \mathbb{R}^n)$.

\section{Main results}\label{Section:Main_results}

In this section we present our main results and their proofs.
The proofs rely on sequences of lemmas, which are proven {in Sections \ref{Section:Proof_of_gauss_conv}--\ref{Section:Proof_of_sup_conv}}.

\subsection{Light traffic.} We first consider the light-traffic case, i.e.,
the case in which $\vk r(u)\to\vk\infty$ as $u\to\infty$. 
Let, {for $\lambda$ defined via {\bf L1},} $\vk B_{\lambda}(\cdot)$ be an $n$-dimensional centered Gaussian process with covariance matrix function $\Sigma(\cdot,\cdot)$ with elements (for $t,s\in\R$, $i,j\in\sprod{n}$)
        \begin{align}
        \Sigma_{ij}(t,s) := \frac{t^{2\lambda} + \left(\mathfrak{r}_{ij}^{-\kappa}s\right)^{2\lambda} - \abs{\,t-\mathfrak{r}_{ij}^{-\kappa}s\,}^{2\lambda}}{2\,\mathfrak{r}^{-\kappa\lambda}_{ij}}\,{\boldsymbol 1}_{\{\mathfrak{r}_{ij}>0\}},\label{sigma_def}
    \end{align}
    where we throughout write $\kappa:=(1-\lambda)^{-1}$.
    Define, for any $t\in[0,\infty)$, the process {$\vk{\mathfrak{X}}_{\lambda}(\cdot)$} by
    \begin{align}
        {\vk{\mathfrak{X}}_{\lambda}}(t)&:=\sup_{-\infty < s < t}\bigl(\vk{\mathcal{C}}\circ\bigl(\vk B_{\lambda}(t) - \vk B_{\lambda}(s)\bigr) - (t-s)\vk 1\bigr)\bigr),\label{X_frac_def}
    \end{align}
    cf.\ \eqref{eq:Ci}.
   Let  $P^{\star}=\{p^\star_{ij}\}_{i,j=1}^n$ be an upper-triangular  matrix with
    \begin{align}
        p^{\star}_{ij} := p_{ij}\,\mathfrak{r}_{ij}^{\kappa\lambda}.\label{P*_def}
    \end{align}

{The following theorem constitutes the main finding of this section.}
\begin{theo} \label{main_light_traffic}
Assume that $\vk r(u)\to \infty$ as $u\to\infty$.
    Let the model satisfy \textbf{N1}, \textbf{N2}, \textbf{N3}, \textbf{L1} and \textbf{L2}.
    Then, 
    \begin{align*}
        \vk Q^{\delta}_u(\cdot) \to {\vk{\mathfrak{Q}}_{\lambda}(\cdot)}:=(I-P^{\star})^{\top}{\vk{\mathfrak{X}}_{\lambda}(\cdot)},
    \end{align*}
    as $u\to\infty$, in distribution in $C([0,\infty)\to\R^n)$.
\end{theo}

\begin{remark}\label{rem:light}\em Using \eqref{sigma_def} we can present the structure of the process $\vk B_{\lambda}(\cdot)$ in a more transparent way. Partition the set $\sprod{n}$ into equivalence classes $\mathcal{I}_1, \ldots, \mathcal{I}_m$ with respect to the asymptotic behavior of $r_i(u)$, meaning that $i,j\in\sprod{n}$ belong to the same class if and only if $\mathfrak{r}_{ij}>0$. For each $i\in\sprod{n}$ we define the index of its class by $f(i)$ (i.e. $i\in\mathcal{I}_{f(i)}$), and for each class $\mathcal{I}_i$ we denote its minimal element by $\mathcal{l}(i)$. For each $i\in\sprod{m}$ we define by $B^{\prime}_{i}(\cdot)$
an independent copy of a fractional Brownian motion with Hurst parameter {$\lambda\in(0,1)$},
{that is a centered Gaussian process with stationary increments and ${\rm Var}\, B_{\lambda}(t)= t^{2\lambda}$ }. Define a new $n$-dimensional centered Gaussian process $\vk B^{\star}(\cdot)$, for $t\in\R$, as follows: with ${\mathcal k}(i):= \mathcal{l}(f(i))$,
    \begin{align*}
        B^{\star}_i(t) := \frac{B^{\prime}_{f(i)}\left(\mathfrak{r}^{-\kappa}_{\mathcal{k}(i),i}t\right)}{\mathfrak{r}^{-\kappa\lambda}_{\mathcal{k}(i),i}}.
    \end{align*}
We now show that the covariance matrix of $\vk B^{\star}(\cdot)$ equals to $\Sigma(\cdot, \cdot)$.
Fix  $i,j\in\langle n\rangle$ such that $i\leqslant j$.
If $\mathfrak{r}_{ij}=0$, then $i$ and $j$ {belong} to different equivalence classes, meaning that $f(i)\not =f(j)$, and hence the processes $B^{\star}_i(\cdot)$ and $B^{\star}_j(\cdot)$ are independent. As a consequence, for such pairs $i,j$ we have, for any $t,s\in\R$,
\begin{align*}
    \operatorname{Cov}\left(B^{\star}_i(t),B^{\star}_j(s)\right)=0.
\end{align*}
Consider now the case $\mathfrak{r}_{ij}>0$. Then $f(j)=f(i)$, so that, for any $t,s\in\R$,
\begin{align*}
    \operatorname{Cov}\left(B^{\star}_i(t),\,B^{\star}_j(s)\right) &= \operatorname{Cov}\left(\frac{B^{\prime}_{f(i)}\left(\mathfrak{r}^{-\kappa}_{\mathcal{k}(i),i}t\right)}{\mathfrak{r}^{-\kappa\lambda}_{\mathcal{k}(i),i}},\,\frac{B^{\prime}_{f(i)}\left(\mathfrak{r}^{-\kappa}_{\mathcal{k}(i),j}s\right)}{\mathfrak{r}^{-\kappa\lambda}_{\mathcal{k}(i),j}}\right)\\
    &=\frac{\operatorname{Cov}\left(B^{\prime}_{f(i)}\left(\mathfrak{r}^{-\kappa}_{\mathcal{k}(i),i}t\right),\,B^{\prime}_{f(i)}\left(\mathfrak{r}^{-\kappa}_{\mathcal{k}(i),j}s\right)\right)}{\mathfrak{r}^{-\kappa\lambda}_{\mathcal{k}(i),i}\mathfrak{r}^{-\kappa\lambda}_{\mathcal{k}(i),j}}\\
    &=\frac{\left(\mathfrak{r}^{-\kappa}_{\mathcal{k}(i),i}t\right)^{2\lambda} + \left(\mathfrak{r}^{-\kappa}_{\mathcal{k}(i),j}s\right)^{2\lambda} - \abs{\mathfrak{r}^{-\kappa}_{\mathcal{k}(i),i}t - \mathfrak{r}^{-\kappa}_{\mathcal{k}(i),j}s}^{2\lambda}}{2\mathfrak{r}^{-\kappa\lambda}_{\mathcal{k}(i),i}\mathfrak{r}^{-\kappa\lambda}_{\mathcal{k}(i),j}}\\
    &=\frac{t^{2\lambda} + \left({\displaystyle \frac{\mathfrak{r}^{-\kappa}_{\mathcal{k}(i),j}}{\mathfrak{r}^{-\kappa}_{\mathcal{k}(i),i}}}s\right)^{2\lambda} - \abs{\,t - {\displaystyle \frac{\mathfrak{r}^{-\kappa}_{\mathcal{k}(i),j}}{\mathfrak{r}^{-\kappa}_{\mathcal{k}(i),i}}}s\,}^{2\lambda}}{2\mathfrak{r}^{-\kappa\lambda}_{\mathcal{k}(i),j} / \mathfrak{r}^{-\kappa}_{\mathcal{k}(i),i}}
    =\frac{t^{2\lambda} + \left(\mathfrak{r}_{ij}^{-\kappa}s\right)^{2\lambda} - \abs{t-\mathfrak{r}_{ij}^{-\kappa}s}^{2\lambda}}{2\mathfrak{r}^{-\kappa\lambda}_{ij}},
\end{align*}
where in the last equality we used that $\mathcal{k}(i)\leqslant i \leqslant j$, implying that $\mathfrak{r}_{\mathcal{k}(i),j} / \mathfrak{r}_{\mathcal{k}(i),i} = \mathfrak{r}_{ij}$. It thus follows that the covariance matrices of $\vk B_{\lambda}(\cdot)$ and $\vk B^{\star}(\cdot)$ coincide, entailing that these two processes are distributionally equivalent. \hfill    ${\Diamond}$
\end{remark}

We proceed by presenting two examples. In the first, all queues decouple, in that they become independent as $u \to \infty$; cf.\ the results presented in \cite{DAI}. In the second example, a subset of the queues remains dependent in the limit.

\begin{example}\label{E1} \em Consider a three-node network, where the first node is fed by the external input
$J(\cdot)$ that satisfies assumptions \textbf{L1}, \textbf{L2}, and the output of the first node is split between the second and third nodes in the proportion $p:1-p$ with $p\in(0,1)$ being fixed. The output rates are `power-law':
\begin{align*}
    r_1(u) = u^{\alpha_1},\qquad r_2(u) = u^{\alpha_2},\qquad {r_3(u) = u^{\alpha_3}}
\end{align*}
for some ${\alpha_1,\alpha_2,\alpha_3>0}$, where we assume ${\alpha_1>\alpha_2>\alpha_3}$. In this case
\begin{align*}
    &\mathcal{C}_1=1,\qquad \mathcal{C}_2 = p,\qquad \mathcal{C}_3 = (1-p),\qquad{\mathfrak{r}_{1,2} = \mathfrak{r}_{1,3} = \mathfrak{r}_{2,3} = 0.}
\end{align*}
Hence, for any $t\geqslant 0$,
\begin{align*}
    {\mathfrak{X}_{\lambda,1}(t)} &= \sup_{-\infty<s<t}\bigl(\bigl(B_{\lambda,1}(t) - B_{\lambda,1}(s)\bigr)- (t - s)\bigr),\\
    {\mathfrak{X}_{\lambda,2}(t)} &= \sup_{-\infty<s<t}\bigl(p\,\bigl(B_{\lambda,2}(t) - B_{\lambda,2}(s)\bigr)- (t-s)\bigr),\\
    {\mathfrak{X}_{\lambda,3}(t)} &= \sup_{-\infty<s<t}\bigl((1-p)\bigl(B_{\lambda,3}(t) -B_{\lambda,3}(s)\bigr) -(t- s)\bigr),
\end{align*}
where $B_{\lambda,1}(\cdot)$, $B_{\lambda,2}(\cdot)$ and $B_{\lambda,3}(\cdot)$ are mutually independent fractional Brownian motions with Hurst parameter $\lambda$. {The  matrix $P^{\star}$, as defined in \eqref{P*_def}, is readily evaluated:
\begin{alignat*}{4}
  &p^{\star}_{11} = p^{\star}_{21} = p^{\star}_{22} = p^{\star}_{31} = p^{\star}_{32} = p^{\star}_{33} = 0,&& && p^{\star}_{12} = p_{12} \mathfrak{r}_{12}^{\kappa\lambda} = p\cdot 0 = 0,\\
  & p^{\star}_{13} = p_{13} \mathfrak{r}_{13}^{\kappa\lambda} = (1-p)\cdot 0 = 0, && \qquad &&p^{\star}_{23} = p_{23} \mathfrak{r}_{23}^{\kappa\lambda} = 0\cdot 0 = 0.
\end{alignat*}
Thus, $P^{\star}$ is an all-zeroes matrix, and hence}, due to \netheo{main_light_traffic},
\begin{align*}
    {\vk{\mathfrak{Q}}_{\lambda}(\cdot) = (I-P^{\star})^{\top}\vk{\mathfrak{X}}_{\lambda}(\cdot) = \vk{\mathfrak{X}}_{\lambda}(\cdot).}
\end{align*}
{Conclude that the three limiting workload processes  are mutually independent; cf.\ the results in \cite{DAI}.}
\hfill $\clubsuit$
\end{example}

\begin{example} \em We consider a three-node network that slightly differs from the one studied in Example~\ref{E1}, but which displays intrinsically different behavior. The only difference is that  the output rates are now
\begin{align*}
    r_1(u) = u^{\alpha_1},\qquad r_2(u) = u^{\alpha_2},\qquad r_3(u) = cu^{\alpha_2}
\end{align*}
for some $c,\alpha_1,\alpha_2>0$ with $\alpha_1>\alpha_2$ and $c\in[0,1]$. Hence, the output rates of the second and third nodes are, up to the constant $c$, asymptotically equivalent. We thus find
\begin{align*}
    &\mathcal{C}_1=1,\qquad \mathcal{C}_2 = p,\qquad \mathcal{C}_3 = (1-p),\qquad\mathfrak{r}_{1,2} = \mathfrak{r}_{1,3} = 0,\qquad \mathfrak{r}_{2,3} = c.
\end{align*}
Hence, for any $t\geqslant 0$, both $\mathfrak{X}_{\lambda,1}(\cdot)$ and $\mathfrak{X}_{\lambda,2}(\cdot)$ are as in Example \ref{E1}, again with $B_{\lambda,1}(\cdot)$ and $B_{\lambda,2}(\cdot)$ being two independent fractional Brownian motions with Hurst parameter $\lambda$. For $\mathfrak{X}_{\lambda,3}(\cdot)$, however, we have
\begin{align*}
    \mathfrak{X}_{\lambda,3}(t) &= \sup_{-\infty<s<t}\left((1-p)\left(\frac{B_{\lambda,2}\left(c^{{-\kappa}}t\right)}{c^{{-\kappa\lambda}}} -\frac{B_{\lambda,2}\left(c^{{-\kappa}}s\right)}{c^{{-\kappa\lambda}}} \right)- (t-s)\right);
\end{align*}
observe in particular that $\mathfrak{X}_{\lambda,2}(\cdot)$ and $\mathfrak{X}_{\lambda,3}(\cdot)$ are
{\it dependent}, as they are driven by the {\it same} process $B_{\lambda,2}(\cdot)$.
In the same way as above, it is seen that $P^{\star}$ is an all-zeroes matrix. Hence,
appealing  to \netheo{main_light_traffic}, we again obtain that
$\vk{\mathfrak{Q}}_{\lambda}(\cdot) = (I-P^{\star})^{\top}\vk{\mathfrak{X}}_{\lambda}(\cdot)
= \vk{\mathfrak{X}}_{\lambda}(\cdot)$, but now there is dependence between the limiting stationary
workloads of the two downstream queues (while there is independence between them and the upstream queue).
\hfill $\clubsuit$
\end{example}

\prooftheo{main_light_traffic} \COM{We start by stating an auxiliary result which provides us with a uniform bound on the function $\mathcal{F}_x(\cdot)$, as defined in \eqref{F_cal_def}.
The proof of this technical lemma can be found in the appendix.}{First} we establish, in Lemma \ref{gauss_convergence} below, a convergence result for the Gaussian processes we take the supremum over in \eqref{X_bar_def}; its proof is postponed to Section \ref{Section:Proof_of_gauss_conv}. Define, for  $t\in\R$,
\begin{align}
    \vk J^{\delta}_{u}(t):=\left(\frac{J(\delta(r_1(u))t)}{\sigma(\delta(r_1(u)))},\ldots, \frac{J(\delta(r_n(u))t)}{\sigma(\delta(r_n(u)))}\right)^{\top}.\label{J_delta_def}
\end{align}

\begin{lem}\label{gauss_convergence}
{Let the model satisfy \textbf{L1}, \textbf{L2},
\textbf{N1}, \textbf{N2} and \textbf{N3}.} Then
    \begin{align}
        \vk J^{\delta}_{u}(\cdot) \to \vk B_{\lambda}(\cdot),\label{eq:gauss_convergence}
    \end{align}
    as $u\to\infty$, in distribution in $C([-T,T]\to\R^n)$, for any $T>0$.
\end{lem}

{The next step is to show, in Lemma \ref{sup_convergence},
the functional convergence of the process 
$\overline{\vk X}_u^{\delta}(\cdot)$, as defined in \eqref{X_delta_def},}
in $C([0,\infty)\to\R^n)$.
{Its proof, as provided in Section \ref{Section:Proof_of_sup_conv},
uses Lemma \ref{gauss_convergence}; it is in this proof that one has to deal with the {subtleties}
arising from the discontinuity of the supremum functional in $C([0,\infty)\to {\mathbb R}^n)$.}
\begin{lem} \label{sup_convergence}
    The following convergence holds in distribution in $C([0,\infty)\to\R^n)$: as $u\to\infty$,
    \begin{align*}
          \overline{\vk X}_u^{\delta}(\cdot)\to {\vk{\mathfrak{X}}_{\lambda}(\cdot)}.
    \end{align*}
\end{lem}
Having Lemma \ref{sup_convergence} at our {disposal}, we can proceed by providing the proof of  \netheo{main_light_traffic}. As a consequence of definitions \eqref{Q_vk_def} and \eqref{X_delta_def}, we can write
\begin{align*}
    \vk Q^{\delta}_u(t) &= \operatorname{diag}\biggl(1/\sigma(\delta(\vk r(u)))\biggr)(I-P^{\top})\operatorname{diag}\biggl(\sigma(\delta(\vk r(u)))\biggr)\operatorname{diag}\biggl(1/\sigma(\delta(\vk r(u)))\biggr) \overline{\vk X}_u(t).
\end{align*}
Hence, applying \nelem{sup_convergence}, to verify the claim of \netheo{main_light_traffic} it suffices to show that, as $u\to\infty$,
\begin{align*}
    \operatorname{diag}\biggl(1/\sigma(\delta(\vk r(u)))\biggr)(I-P^{\top})\operatorname{diag}\biggl(\sigma(\delta(\vk r(u)))\biggr) \to (I-P^{\star})^{\top},
\end{align*}
where the matrix $P^*$ is  defined in \eqref{P*_def}. This means that, due to \eqref{delta_def}, it is enough to show that, for any $i,j\in\sprod{n}$,
\begin{align}
     \lim_{u\to\infty} P_{ij}\frac{r_i(u)\,\delta(r_i(u))}{r_j(u)\,\delta(r_j(u))} = P^{\star}_{ij}.\label{P_bar_def}
\end{align}
Using that the  matrix $P$ is strictly upper-triangular,
we can restrict ourselves to $i<j$.\\
{\bf {Case 1:}  $\mathfrak{r}_{ij}>0$.}
Using that $\delta(x)\in\operatorname{\mathscr{RV}}_{\infty}({-\kappa})$ with  $\kappa=(1-\lambda)^{-1}$,
{we have}
\begin{align*}
 \lim_{u\to\infty} P_{ij}\frac{\sigma\bigl(\delta(r_i(u))\bigr)}{\sigma\bigl(\delta(r_j(u))\bigr)} = \lim_{u\to\infty}P_{ij}\frac{r_i(u)\,\delta(r_i(u))}{r_j(u)\,\delta(r_j(u))} = \mathfrak{r}_{ij}^{-1}P_{ij}\lim_{u\to\infty} \frac{\delta\left(\frac{r_i(u)}{r_j(u)}\,r_j(u)\right)}{\delta(r_j(u))} = \mathfrak{r}_{ij}^{{\kappa-1}}P_{ij} =  P_{ij}^{\star}.
\end{align*}
{\bf {Case 2:} $\mathfrak{r}_{ij}=0$.}
Denote, for $i,j\in\sprod{n}$ and $u>0$,
\begin{align}
    \Delta_{ij}(u) := \delta(r_i(u))/\delta(r_j(u)).\label{Delta_ij_def}
\end{align} In view of \eqref{delta_rv}, one can derive (for instance, by applying the Karamata representation theorem \cite[(0.35)]{resnick2008extreme}) that
$\Delta_{ij}(u) \to 0$ as $u\to\infty$.
Thus,
\begin{align*}
    \lim_{u\to\infty} P_{ij}\frac{\sigma\bigl(\delta(r_j(u))\bigr)}{\sigma\bigl(\delta(r_i(u))\bigr)} = 0 = P^{\star}_{ij}.
\end{align*}
Hence \eqref{P_bar_def} follows, {which completes} the proof of \netheo{main_light_traffic}.
\QED

\subsection{Heavy traffic} This subsection covers the heavy-traffic case, i.e.,
the case in which $\vk r(u)\to\vk 0$ as $u\to\infty$. We start by introducing some notation needed to state the main result; these objects align with their light-traffic counterparts.
Let, {for $\alpha$ defined in {\bf H1},} $\vk B_{\alpha}(\cdot)$ denote an $n$-dimensional centered Gaussian process with covariance matrix function $\Sigma(\cdot,\cdot)$ with elements (for $t,s\in\R$,  $i,j\in\sprod{n}$)
    \begin{align*}
        \Sigma_{ij}(t,s) := \frac{t^{2\alpha} + \left(\mathfrak{r}_{ij}^{\xi}s\right)^{2\alpha} - \abs{t-\mathfrak{r}_{ij}^{\xi}s}^{2\alpha}}{2\mathfrak{r}^{\xi \alpha}_{ij}}\,{\boldsymbol 1}_{\{\mathfrak{r}_{ij}>0\}},
    \end{align*}
    where $\xi:=1/(\alpha-1)$. For $t>0$, {we define}
    \begin{align*}
        {\vk{\mathfrak{X}}_{\alpha}(t)}&:=\sup_{-\infty < s < t}\bigl(\vk{\mathcal{C}}\circ\bigl(\vk B_{\alpha}(t) -\vk B_{\alpha}(s)\bigr)- (t-s)\vk 1\bigr) .
    \end{align*}
    Let  $P^{\star}=\{p^\star_{ij}\}_{i,j=1}^n$ be an upper-triangular  matrix with
    \begin{align*}
        p^{\star}_{ij} := p_{ij}\,\mathfrak{r}_{ij}^{\xi\alpha}.
    \end{align*}

\begin{theo} \label{main_heavy_traffic}
    Assume that $\vk r(u)\to 0$ as $u\to\infty$.
    Let the model satisfy \textbf{N1}, \textbf{N2}, \textbf{N3}, \textbf{H1} and \textbf{H2}.
    Then,
    \begin{align*}
        \vk Q^{\delta}_u(\cdot) \to {\vk{\mathfrak{Q}}_{\alpha}(\cdot)}:=(I-P^{\star})^{\top}{\vk{\mathfrak{X}}_{\alpha}(\cdot)},
    \end{align*}
as $u\to\infty$, in distribution in $C([0,\infty)\to\R^n)$.
    \end{theo}

\begin{remark}\label{rem:heavy}\em {Similarly to} Remark \ref{rem:light}, the process $\vk B_{\alpha}(\cdot)$ can be represented in an alternative form. Indeed, we can write, for $t\in\R$,
\begin{align*}
        B^\star_{i}(t) = \frac{B^{\prime}_{f(i)}\left(\mathfrak{r}^{\xi}_{\mathcal{k}(i),i}t\right)}{\mathfrak{r}^{\xi \alpha}_{\mathcal{k}(i),i}},
\end{align*}
 with, for each $i\in\sprod{m}$, $B^{\prime}_{i}(\cdot)$ now denoting an independent copy of a fractional Brownian motion with Hurst parameter $\alpha$.
 As was  observed in Remark \ref{rem:light}, one can argue that $\vk B_{\alpha}(\cdot)$ and $\vk B^\star(\cdot)$ are distributionally equivalent. \hfill    ${\Diamond}$
\end{remark}

\prooftheo{main_heavy_traffic}
The proof of this theorem follows  by using an argumentation that closely resembles the one used
in the proof of \netheo{main_light_traffic}.
Let {$\vk J_{u}^{\delta}(\cdot)$} be as defined in \eqref{J_delta_def}.
The counterpart of \nelem{gauss_convergence} can be stated as follows.
Its proof is fully identical to the proof of \nelem{gauss_convergence} and therefore is omitted.

\begin{lem}\label{gauss_convergence_heavy}
{Let the model satisfy \textbf{H1}, \textbf{H2},
\textbf{N1}, \textbf{N2} and \textbf{N3}}. Then
    \begin{align*}
        \vk J^{\delta}_{u}(\cdot) \to \vk B_{\alpha}(\cdot),
    \end{align*}
as $u\to\infty$, in distribution in $C([-T,T]\to\R^n)$, for any $T>0$.
\end{lem}

The next step is the counterpart of \nelem{sup_convergence},
whose proof follows by the same steps as the proof of \eqref{sup_convergence}.

\begin{lem} \label{sup_convergence_heavy}
    As $u\to\infty$ the following convergence holds in distribution in $C([0,\infty]\to\R^n)$
    \begin{align*}
          \overline{\vk X}_u^{\delta}(\cdot)\to {\vk{\mathfrak{X}}_{\alpha}(\cdot)},
    \end{align*}
    where $ \overline{\vk X}_u^{\delta}(\cdot)$ is defined in \eqref{X_delta_def} and ${\vk{\mathfrak{X}}_{\alpha}(\cdot)}$ is defined in \netheo{main_heavy_traffic}.
\end{lem}

The rest of the proof of Theorem \ref{main_heavy_traffic} is analogous to the proof of \netheo{main_light_traffic}.
\QED
\begin{remark}\label{rem_deg} \em
We remark that
although the proofs of Theorems \ref{main_light_traffic} and \ref{main_heavy_traffic}
follow by highly similar arguments and the classes of the limiting random structures coincide,
there is a subtle but important difference.
As noted in Remark \ref{rem:H2_L2}, the heavy-traffic case allows for
driving processes with standard deviation functions that grows linearly at $0$,
while the light-traffic case excludes such processes.
The class of integrated stationary Gaussian processes contains natural examples of such processes; see e.g. \cite{DeM04}.
For an explanation of this symmetry breaking, we consider a single queue ($n=1$), and we
suppose that $J(t)=\int_0^t \eta(s)\,{\rm d}s$, where $\eta(\cdot)$ is a stationary centered Gaussian process
with continuous sample paths a.s.\ 
\begin{itemize}
    \item[$\circ$] We start by analyzing the light-traffic scenario
($r_1(u)\to \infty$ as $u\to\infty$).
Observe that {\bf L2} is satisfied.  In the setting considered the standard deviation function of $J(\cdot)$
is such that $\sigma(t)\sim Rt$ as $t\downarrow 0$, with $R\in (0,\infty)$.
In order to obtain an analog of
Lemma \ref{gauss_convergence}, suppose that $\delta(x)\to 0$, as $x\to \infty$.
By the same argument as in the proof of Lemma \ref{gauss_convergence},
see also the proof of Theorem 4.2 in \cite{DeM04},
$J_u^{\delta}(\cdot)\to B_1(\cdot)$, as $u\to\infty$.
Unfortunately, this convergence does not extend to the convergence of
\[
X_{u,1}^{\delta}(0)=_{\rm d}\sup_{s\geqslant 0}\left(J_u^{\delta}(s)- \frac{r_1(u)\delta(r_1(u))}{\sigma(\delta(r_1(u)))}s\right),
\]
since 
\[
\frac{r_1(u)\delta(r_1(u))}{\sigma(\delta(r_1(u)))}\sim
\frac{r_1(u)\delta(r_1(u))}{R\delta(r_1(u))}\sim \frac{r_1(u)}{R}\to \infty,
\]
as $u \to \infty$, which leads to a degenerate limit for $X_{u,1}^{\delta}(0)$ in the regime that $u\to\infty$.
\item[$\circ$]
In the heavy-traffic scenario ($r_1(u)\to 0$ as $u\to\infty$), the situation is different, since the limit process depends on $\alpha$ in
{\bf H1}, which is responsible for the asymptotic behavior of
$\sigma(\cdot)$ at $\infty$. Because {\bf H2} is satisfied with $\beta=1$ in this case,
Theorem \ref{main_heavy_traffic} can be applied to
the family of integrated stationary Gaussian processes as long as $\sigma(\cdot)$ satisfies {\bf H1}.
\end{itemize}
We refer to \cite{DeM04, Gaussianqueues} for related models where the observations also apply. 
\hfill$\Diamond$
\end{remark}

\section{Proof of Lemma \ref{gauss_convergence}}\label{Section:Proof_of_gauss_conv}

To establish the claim of Lemma \ref{gauss_convergence}, i.e.,
the convergence in distribution \eqref{eq:gauss_convergence},
we rely on the approach described in \cite[Section 7]{billingsley2013convergence}.
It consists of {(i)~the proof of the convergence of the finite-dimensional distributions and (ii)~the verification of the tightness condition.}
More precisely, in the first step we are to show that, for any
$-T\leqslant t_1<\ldots<t_m\leqslant T$, as $u\to\infty$,
\begin{align}
    \left((\vk J_u^{\delta}(t_1))^{\top},\ldots, (\vk J_u^{\delta}(t_m))^{\top}\right)^{\top}\to \left((\vk B_{\lambda}(t_1))^{\top},\ldots, (\vk B_{\lambda}(t_m))^{\top}\right)^{\top}\label{fdd_for_J}
\end{align}
in distribution. In the second step we prove that the sequence $\vk J_u^{\delta}(\cdot)$ is tight in $C([-T,T]\to\R^n)$. Due to \cite[Theorem 7.5]{billingsley2013convergence}, {given \eqref{fdd_for_J},} it suffices to verify that for, any $\eta>0$,
\begin{align}
    \lim_{\zeta \downarrow 0}\limsup_{u\to\infty}\pk{\sup_{\substack{t,x\in[-T,T], \\ \abs{t-s}\leqslant\zeta}}\lvert\lvert\vk J^{\delta}_{u}(t) - \vk J^{\delta}_{u}(s)\rvert\rvert\geqslant\eta}=0.\label{uniformly_equicontinuous}
\end{align}

\subsection{Finite-dimensional distributions} In this subsection we verify \eqref{fdd_for_J}. {First suppose that $t_i=0$ for some $i\in\sprod{m}$.
Observe that \textbf{L1} gives us $\sigma(0)=0$, i.e., $J(0)=0$ almost surely.
It now follows, from~\eqref{J_delta_def},
that $\vk J_u^{\delta}({0}) = \vk 0$ almost surely.
Hence, as $u\to\infty$, $\vk J_u^{\delta}({0})\to \vk B_{\lambda}(0)$ almost surely. Conclude that in~\eqref{fdd_for_J} we can disregard the elements $\vk J_u^{\delta}(t_i)$ for which $t_i=0$ and consider only $t_1,\ldots,t_m\not= 0$.}

Using that the random variables on the left-hand side of \eqref{fdd_for_J} are centered Gaussian for any $u>0$, the limiting vector is Gaussian as well. Hence, to show the convergence in \eqref{fdd_for_J}, it suffices to verify that the covariance matrices of the vectors on the left-hand side converge to the covariance matrix of the vector on the right-hand side.  Fix $i,j\in\langle n\rangle$ such that $i\leqslant j$, and $t,s\in[-T,T]\setminus\{0\}$. Then, using that $J(\cdot)$ has stationary increments,
\begin{align}
    \operatorname{Cov}\left(J^{\delta}_{u,i}(t), J^{\delta}_{u,j}(s)\right) &= \operatorname{Cov}\left(\frac{J\bigl(\delta(r_i(u))t\bigr)}{\sigma\bigl(\delta(r_i(u))\bigr)}, \frac{J\bigl(\delta(r_j(u))s\bigr)}{\sigma\bigl(\delta(r_j(u))\bigr)}\right)\notag\\
    &= \frac{\sigma^2\bigl(\delta(r_i(u))t\bigr)+\sigma^2\bigl(\delta(r_j(u))s\bigr)-\sigma^2\bigl(\abs{\delta(r_i(u))t-\delta(r_j(u))s}\bigr)}{2\sigma\bigl(\delta(r_i(u))\bigr)\sigma\bigl(\delta(r_j(u))\bigr)}\label{Cov_init_representation}
\end{align}
We consider two cases separately: either $\mathfrak{r}_{ij}>0$ or $\mathfrak{r}_{ij}=0$.

\textbf{Case 1: $\mathfrak{r}_{ij}>0$}.
In this case we know that $r_j(u)/r_{i}(u)$ is bounded away from zero for sufficiently large $u$.
Hence, using the definition of $\Delta_{ji}(\cdot)$ given in \eqref{Delta_ij_def}, and by virtue of  \eqref{delta_rv},
\begin{align*}
    \lim_{u\to\infty}{\Delta_{ji}(u)} = \lim_{u\to\infty}\frac{\delta\left({\displaystyle\frac{r_{j}(u)}{r_{i}(u)}}\,r_{i}(u)\right)}{\delta(r_{i}(u))} = \mathfrak{r}^{{-\kappa}}_{ij}>0.
\end{align*}
Hence, applying \eqref{delta_def} we obtain that
\begin{align}
    \lim_{u\to\infty}\frac{\sigma(\delta(r_j(u)))}{\sigma(\delta(r_i(u)))} = \lim_{u\to\infty}\frac{r_j(u)\,\delta(r_j(u))}{r_i(u)\,\delta(r_i(u))} = \mathfrak{r}^{{-\kappa\lambda}}_{ij}>0.\label{sigma_frac_r>0}
\end{align}
Applying \eqref{sigma_frac_r>0} and \eqref{RV_statement} separately for each term of \eqref{Cov_init_representation}, we conclude that
\begin{align*}
    \lim_{u\to\infty}\frac{\sigma^2\bigl(\delta(r_i(u))t\bigr)}{2\sigma\bigl(\delta(r_i(u))\bigr)\sigma\bigl(\delta(r_j(u))\bigr)} &= \lim_{u\to\infty}\frac{\sigma^2\bigl(\delta(r_i(u))\,t\bigr)}{2\sigma^2\bigl(\delta(r_i(u))\bigr)}\frac{\sigma\bigl(\delta(r_i(u))\bigr)}{\sigma\bigl(\delta(r_j(u))\bigr)} \\
    &= \lim_{u\to\infty}\frac{\mathcal{F}^2_{\delta(r_i(u))}(t)}{2\mathfrak{r}^{{-\kappa\lambda}}_{ij}} = \frac{t^{2\lambda}}{2\mathfrak{r}^{{-\kappa\lambda}}_{ij}},\\
    \lim_{u\to\infty}\frac{\sigma^2\bigl(\delta(r_j(u))\,s\bigr)}{2\sigma\bigl(\delta(r_i(u))\bigr)\sigma\bigl(\delta(r_j(u))\bigr)} &= \lim_{u\to\infty}\frac{\sigma^2\bigl(\delta(r_j(u))s\bigr)}{2\sigma^2\bigl(\delta(r_j(u))\bigr)}\frac{\sigma\bigl(\delta(r_j(u))\bigr)}{\sigma\bigl(\delta(r_i(u))\bigr)} \\
    &= \lim_{u\to\infty}\frac{\mathcal{F}^2_{\delta(r_j(u))}(s)\mathfrak{r}^{{-\kappa\lambda}}_{ij}}{2} = \frac{s^{2\lambda}\mathfrak{r}^{{-\kappa\lambda}}_{ij}}{2},\\
    \lim_{u\to\infty}\frac{\sigma^2\bigl(\abs{\delta(r_i(u))t-\delta(r_j(u))s}\bigr)}{2\sigma\bigl(\delta(r_i(u))\bigr)\sigma\bigl(\delta(r_j(u))\bigr)} &= \lim_{u\to\infty}\frac{\sigma^2\bigl(\abs{\delta(r_i(u))t-\delta(r_j(u))s}\bigr)}{2\sigma^2\bigl(\delta(r_i(u))\bigr)}\frac{\sigma\bigl(\delta(r_i(u))\bigr)}{\sigma\bigl(\delta(r_j(u))\bigr)}\\
    &= \lim_{u\to\infty}\frac{\mathcal{F}^2_{\delta(r_i(u))}\left(\abs{t-{\displaystyle \Delta_{ij}(u)
    }s}\right)}{2\mathfrak{r}^{{-\kappa\lambda}}_{ij}} = \frac{\abs{t-\mathfrak{r}_{ij}^{{-\kappa}}s}^{2\lambda}}{2\mathfrak{r}^{{-\kappa\lambda}}_{ij}}.
\end{align*}
Thus, according to \eqref{sigma_def}, the convergence in distribution \eqref{fdd_for_J} follows, in the case that $\mathfrak{r}_{ij}>0$.

\textbf{Case 2: $\mathfrak{r}_{ij}=0$}. In this case, appealing  to \eqref{delta_rv}, we obtain, {using for instance the Karamata representation theorem  \cite[(0.35)]{resnick2008extreme}},
\begin{align}
    \lim_{u\to\infty}\Delta_{ij}(u)= 0. \label{delta_frac}
\end{align}
Consequently, by  \textbf{L1}, along the same lines {as in (\ref{sigma_frac_r>0})},
\begin{align}
    \lim_{u\to\infty}\frac{\sigma\bigl(\delta(r_j(u))\bigr)}{\sigma\bigl(\delta(r_i(u))\bigr)}= 0. \label{sigma_frac_r=0}
\end{align}
The expression for the covariance, as displayed in \eqref{Cov_init_representation}, can be rewritten as ${\mathscr G}_{ij}(u)\,{\mathscr H}_{ij}(u)+{\mathscr I}_{ij}(u)$, where
\begin{align}
    {\mathscr G}_{ij}(u)&:= \frac{\sigma\bigl(\delta(r_i(u))t\bigr)-\sigma\bigl(\abs{\delta(r_i(u))t-\delta(r_j(u))s}\bigr)}{\sigma\bigl(\delta(r_j(u))\bigr)},\label{Cov_1}\\
    {\mathscr H}_{ij}(u)&:= \left(\frac{1}{2}\mathcal{F}_{\delta(r_i)}(t) + \frac{\sigma\bigl(\abs{\delta(r_i(u))t-\delta(r_j(u))s}\bigr)}{2\sigma\bigl(\delta(r_i(u))\bigr)\bigr)}\right),\label{Cov_2}\\
    {\mathscr I}_{ij}(u)&:= \mathcal{F}^2_{\delta(r_j(u))}(s)\frac{\sigma\bigl(\delta(r_j(u))\bigr)}{2\sigma\bigl(\delta(r_i(u))\bigr)}\label{Cov_3}.
\end{align}
We proceed by separately analyzing ${\mathscr G}_{ij}(u)$, ${\mathscr H}_{ij}(u)$, and ${\mathscr I}_{ij}(u)$. Regarding \eqref{Cov_3}, due to \eqref{sigma_frac_r=0},
\begin{align*}
    \limsup_{u\to\infty} {\mathscr I}_{ij}(u)= \limsup_{u\to\infty}\abs{\mathcal{F}^2_{\delta(r_j)}(s)\frac{\sigma\bigl(\delta(r_j(u))\bigr)}{2\sigma\bigl(\delta(r_i(u))\bigr)}} = 0.
\end{align*}
Now consider the factor \eqref{Cov_2}. {Using that $t\not=0$,} for any $\varepsilon\in(0,\abs{t})$, applying \eqref{F_RV_infty} and \eqref{delta_frac} we can derive the following upper bound:
\begin{align*}
\limsup_{u\to\infty} {\mathscr H}_{ij}(u)&= \limsup_{u\to\infty}\abs{\frac{1}{2}\mathcal{F}_{\delta(r_i)}(t) + \frac{\sigma\bigl(\abs{\delta(r_i(u))t-\delta(r_j(u))s}\bigr)}{2\sigma\bigl(\delta(r_i(u))\bigr)\bigr)}}\\ &\leqslant \frac{\abs{t}^{\lambda}}{2} + \frac{1}{2}\limsup_{u\to\infty}\max\limits_{t^{\star}\in [t-\varepsilon,t+\varepsilon]}\mathcal{F}_{\delta(r_i(u))}(\abs{t^*})\\
&= \frac{\abs{t}^{\lambda}}{2} + \max_{t^{\star}\in[t-\varepsilon,t+\varepsilon]}\frac{\abs{t^{\star}}^{\lambda}}{2}=\frac{\abs{t}^{\lambda}}{2} + \frac{\abs{t+\varepsilon}^{\lambda}}{2}.
\end{align*}
Letting  $\varepsilon\downarrow 0$, we obtain that
\[
\limsup_{u\to\infty}{\mathscr H}_{ij}(u)\leqslant \abs{t}^{\lambda}.\]
Finally, for the factor \eqref{Cov_1} we have, {using the compact notation $\Delta^\sigma_{ij}(u) := \sigma(\delta(r_i(u)))/\sigma(\delta(r_j(u)))$,}
\begin{align*}
    \limsup_{u\to\infty}{\mathscr G}_{ij}(u)&=\limsup_{u\to\infty}\frac{\sigma\bigl(\delta(r_i(u))t\bigr)-\sigma\bigl(\delta(r_i(u))t-\delta(r_j(u))s\bigr)}{\sigma\bigl(\delta(r_j(u))\bigr)}\\
    &=\limsup_{{u\to\infty}}\Delta^\sigma_{ij}(u)\left(\mathcal{F}_{\delta(r_i(u))}(t) - \mathcal{F}_{\delta(r_i(u))}\left(t-s\Delta_{ij}(u)\right)\right)\\
    &=\limsup_{{u\to\infty}}\Delta^\sigma_{ij}(u)\,s\,\Delta_{ij}(u)\left(\frac{\mathcal{F}_{\delta(r_i(u))}(t) - \mathcal{F}_{\delta(r_i(u))}\left(t-s{\displaystyle \Delta_{ij}(u)}\right)}{s{\displaystyle \Delta_{ij}(u)}}\right)\\&={\operatorname{sign}(t)}\,s\lambda\,\mathfrak{r}_{ij}{\abs{t}}^{\lambda - 1}=0,
\end{align*}
where in the last equality we have used {that $t\not= 0$ (where $\lambda-1<0$) and} \eqref{delta_def}, in combination with {\eqref{RV_statement}}.
Hence, due to \eqref{sigma_def}, \eqref{fdd_for_J} now also follows in the case that $\mathfrak{r}_{ij}=0$.

\subsection{Tightness} {The aim} of this subsection is to verify \eqref{uniformly_equicontinuous}. {This section relies heavily on uniform power-law bounds for the function $\mathcal{F}_x(\cdot)$ which are stated in the following lemma, with its proof provided in the appendix.}
\begin{lem}\label{lem:F_uniform_bound}
    Under \textbf{L1} and \textbf{L2}, there exist constants $\gamma_{0}\in(0,\infty)$, $\gamma_{\infty}\in(0,1)$, $\mathfrak{C}>0$ and $\mathfrak{a}>0$ such that, for any $x\in(0,\mathfrak{a}]$ and  $t>0$,
    \begin{align}
        \mathcal{F}_x(t)\leqslant \mathfrak{C}(t^{\gamma_{\infty}}+t^{\gamma_{0}}).\label{F_uniform_bound}
    \end{align}
    In addition, for any $t>1$,
    \begin{align}
        \mathcal{F}_x(t)\leqslant 2\mathfrak{C}t^{\gamma_{\infty}}.\label{F_uniform_bound_t>1}
    \end{align}
\end{lem}

\begin{remark}\label{rem:F_uniform_bound_heavy} {\em
    The bounds presented in \eqref{F_uniform_bound} and \eqref{F_uniform_bound_t>1} hold as well under the assumptions \textbf{H1} and \textbf{H2} instead of \textbf{L1} and \textbf{L2}, for all $x\in[\mathfrak{a},\infty)$. We provide the detailed argumentation in the appendix.}
\end{remark}

{In order to show \eqref{uniformly_equicontinuous}, we begin by}  noticing that, for any $\zeta,\eta,u>0$,
\begin{align*}
    \pk{\sup_{\substack{t,x\in[-T,T], \\ \abs{t-s}\leqslant\zeta}}\lvert\lvert\vk J^{\delta}_{u}(t) - \vk J^{\delta}_{u}(s)\rvert\rvert\geqslant\eta}&\leqslant \pk{\sup_{\substack{t,x\in[-T,T], \\ \abs{t-s}\leqslant\zeta}}\sum_{i=1}^{n}\abs{\frac{J\bigl(\delta(r_i(u)t)\bigr)}{\sigma\bigl(\delta(r_i(u))\bigr)} - \frac{J\bigl(\delta(r_i(u)s)\bigr)}{\sigma\bigl(\delta(r_i(u))\bigr)}}\geqslant\eta}\\
    &\leqslant \sum_{i=1}^{n}\pk{\sup_{\substack{t,x\in[-T,T], \\ \abs{t-s}\leqslant\zeta}}\abs{\frac{J\bigl(\delta(r_i(u)t)\bigr)}{\sigma\bigl(\delta(r_i(u))\bigr)} - \frac{J\bigl(\delta(r_i(u)s)\bigr)}{\sigma\bigl(\delta(r_i(u))\bigr)}}\geqslant\eta}.
\end{align*}
Hence, to establish \eqref{uniformly_equicontinuous}, it is sufficient to show that, for any $\eta>0$,
\begin{align}
    \lim_{\zeta\downarrow 0}\limsup_{u\to\infty}\pk{{\sup_{(s,t)\in {\mathscr L}_{T,\zeta}}\abs{J^{\star}_{u}(s,t)}\geqslant\eta}} = 0\label{uniformly_equicontinuous_component}
\end{align}
{where we denote
\begin{align*}
    J^{\star}_{u}(s,t) := \frac{J\bigl(\delta(u)t\bigr)}{\sigma\bigl(\delta(u)\bigr)} - \frac{J\bigl(\delta(u)s\bigr)}{\sigma\bigl(\delta(u)\bigr)}
\end{align*}
and ${\mathscr L}_{T,\zeta}:=\{(s,t)\in[T,T]: |t-s|\leqslant \zeta\}.$}
{We recall that we consider the regime that $r_i(u)\to\infty$ as $u\to\infty$}.
Applying \eqref{F_uniform_bound} from \nelem{lem:F_uniform_bound}, it follows that, taking  $u$ sufficiently large to make sure that $\delta(u)<\mathfrak{a}$, for  any {$s,t\in[-T,T]$},
\begin{align}
    \E{{\bigl(J^{\star}_{u}(s,t)\bigr)^2}} &= \mathcal{F}^2_{\delta(u)}(\abs{{t-s}})\leqslant {\mathfrak{C}^2}\left(\abs{{t-s}}^{\gamma_{\infty}} + \abs{{t-s}}^{\gamma_{0}}\right)^2\notag\\
    &= {\mathfrak{C}^2\abs{{t-s}}^{2\min(\gamma_{\infty},\gamma_{0})}\left(1+\abs{{t-s}}^{\max(\gamma_{\infty},\gamma_0)-\min(\gamma_{\infty},\gamma_0)}\right)^{2}}\notag\\
    &\leqslant {\mathfrak{C}^2\abs{{t-s}}^{2\min(\gamma_{\infty},\gamma_{0})}\left(1+(2T)^{\max(\gamma_{\infty},\gamma_0)-\min(\gamma_{\infty},\gamma_0)}\right)^{2}},\notag
\end{align}
{in the second equality applying the obvious identity $x^A+x^B= x^{\min(A,B)}+x^{\max(A,B)}$ (valid for any $x\geqslant 0$ and $A,B\in{\mathbb R}$).} {It implies that for any $s_1,t_1,s_2,t_2\in[-T,T]$, denoting
$
    \mathfrak{C}_2 = \mathfrak{C}\left(1+(2T)^{\max(\gamma_{\infty},\gamma_0)-\min(\gamma_{\infty},\gamma_0)}\right)
$
we obtain for sufficiently large $u$
\begin{align}
    \E{\bigl(J^{\star}_{u}(s_1,t_1) - J^{\star}_u(s_2,t_2)\bigr)^2}&\leqslant 2 \left(\E{\bigl(J^{\star}_{u}(s_1,s_2)\bigr)^2} + \E{\bigl(J^{\star}_{u}(t_1,t_2)\bigr)^2}\right)\notag\\
    &\leqslant 2\mathfrak{C}_2^2\left(\abs{s_2-s_1}^{2\min(\gamma_{\infty},\gamma_{0})} + \abs{t_2-t_1}^{2\min(\gamma_{\infty},\gamma_{0})}\right)\notag\\
    &\leqslant 2^{2-\min(\gamma_{\infty},\gamma_0)}\mathfrak{C}_2^2\left(\abs{s_2-s_1}^2 + \abs{t_2-t_1}^2\right)^{\min(\gamma_{\infty},\gamma_0)},\label{Pit_claim_J}
\end{align}
where in the last line we apply the generalized mean inequality.} Inequality \eqref{Pit_claim_J} allows us to apply Piterbarg's inequality \cite[Theorem 8.1]{piterbarg1996asymptotic}, {which provides us} the following upper bound:
for some constant $C>0$, a sufficiently large $u>0$ and any $\zeta>0$,
\begin{align}
    \pk{\sup_{{(s,t)\in {\mathscr L}_{T,\zeta}}}{\abs{J^{\star}_{u}(s,t)}}\geqslant\eta}\leqslant \varphi(\eta)\,\sigma^{\star}_{u}(\zeta)\,e^{-\frac{\eta^2}{2\left(\sigma^{\star}_{u}(\zeta)\right)^2}}, \label{pit_inequality}
\end{align}
where
\begin{align*}
\varphi(\eta):= {\frac{C(2T)^2\eta^{4/(2\min(\gamma_{\infty},\gamma_0))}}{\sqrt{2\pi}\eta},}\qquad
    \bigl(\sigma^{\star}_{u}(\zeta)\bigr)^2 := \sup_{{(s,t)\in {\mathscr L}_{T,\zeta}}}\E{{\bigl(J^{\star}_{u}(s,t)\bigr)^2}} {\leqslant} \sup_{t\in[0,\zeta]}\mathcal{F}^2_{\delta(u)}(t).
\end{align*}
{By applying \eqref{F_uniform_bound} we obtain, for sufficiently large $u$ such that $\delta(u)<\mathfrak{a}$, where $\mathfrak{a}$ is defined in \nelem{lem:F_uniform_bound},}
\begin{align*}
    {\bigl(\sigma^{\star}_{u}(\zeta)\bigr)^2 \leqslant \sup_{t\in[0,\zeta]}\mathfrak{C}^2\left(t^{\gamma_{\infty}} + t^{\gamma_{0}}\right)^2 = \mathfrak{C}^2(\zeta^{\gamma_{\infty}} + \zeta^{\gamma_{0}})^2.}
\end{align*}
We conclude that $\bigl(\sigma^{\star}_{u}(\zeta)\bigr)^2$ is bounded for large $u$, in that
\begin{align}
    \limsup_{u\to\infty}\,\bigl(\sigma^{\star}_{u}(\zeta)\bigr)^2 \leqslant {\mathfrak{C}^2(\zeta^{\gamma_{\infty}} + \zeta^{\gamma_{0}})^2}.\label{sigma_stal_limit}
\end{align}
Upon combining \eqref{pit_inequality} with \eqref{sigma_stal_limit}, we obtain the following upper bound on \eqref{uniformly_equicontinuous_component}: for any $\eta>0$,
\begin{align*}
    \lim_{\zeta\downarrow 0}\limsup_{u\to\infty}&\,\pk{\sup_{{(s,t)\in {\mathscr L}_{T,\zeta}}}\abs{{J^{\star}_{u}(s,t)}}\geqslant\eta}\\& \leqslant {\varphi(\eta)}\,{\mathfrak{C}}\,\lim_{\zeta\downarrow 0}\left({(\zeta^{\gamma_{\infty}} + \zeta^{\gamma_{0}})}\,\exp\left({-\frac{\eta^2}{{2\,\mathfrak{C}^2(\zeta^{\gamma_{\infty}} + \zeta^{\gamma_{0}})^2}}}\right) \right)=0.
\end{align*}
Hence, \eqref{uniformly_equicontinuous_component} follows.
We finish this proof by remarking that one could alternatively show \eqref{uniformly_equicontinuous_component} essentially following  the line of the proof of \cite[Proposition 4]{Gaussianqueues}. \hfill
\QED

\section{Proof of Lemma \ref{sup_convergence}}\label{Section:Proof_of_sup_conv}
{Lemma \ref{gauss_convergence} established a functional limit theorem of a process that is directly
related to the network's Gaussian input process $J(\cdot)$, i.e., the process-level convergence  $\vk J^{\delta}_{u}(\cdot)\to\vk B_{\lambda}(\cdot)$. In this section we establish Lemma
\ref{sup_convergence}, i.e., the process-level convergence  of the `reflected versions':
$\overline{\vk X}_u^{\delta}(\cdot)\to{\vk{\mathfrak{X}}_{\lambda}(\cdot)}.$
It is noted that the proof of Lemma \ref{sup_convergence} is intrinsically harder than that of Lemma
\ref{gauss_convergence}: bearing in mind the definition of $\overline{\vk X}_u^{\delta}(\cdot)$ in terms of a supremum,
there is the technical challenge of dealing with the discontinuity of the supremum functional in
{$C((-\infty,T]\to {\mathbb R}^n)$.} As in the proof of Lemma \ref{gauss_convergence}, we first
show convergence of finite-dimensional distributions, and then we establish tightness.}

\subsection{Finite-dimensional distributions} Let $m\in\N$. Fix some $t_1,\ldots,t_m\in[0,\infty)$ and $\vk v_1,\ldots,\vk v_m\in\R^{n}$. The objective of this subsection is to show that, as $u\to\infty$,
\begin{align}
    \pk{\begin{matrix}\overline{\vk X}_u^{\delta}(t_1)  >  \vk v_1 \\  \ldots \\ \overline{\vk X}_u^{\delta}(t_m)  >  \vk v_m \end{matrix}}\to \pk{ \begin{matrix}{\vk{\mathfrak{X}}_{\lambda}(t_1)} > \vk v_1 \\ \ldots \\ {\vk{\mathfrak{X}}_{\lambda}(t_m)} > \vk v_m \end{matrix}}.\label{sup_fdd_claim}
\end{align}
The main idea now is, in order to apply \nelem{gauss_convergence}, to decompose the supremum over the interval $(-\infty,t]$, which appears in \eqref{X_bar_def}, into the suprema over two intervals, namely $(-\infty,-T_{\varepsilon})$ and $[-T_{\varepsilon}, t]$. When choosing $T_{\varepsilon}$ large enough, the first supremum is negligible, while for the second supremum   we obtain the convergence \eqref{sup_fdd_claim} by applying the continuous mapping theorem.
So as to show that one can ignore the supremum over $(-\infty,T_{\varepsilon}]$, we state the following lemma, proven in the appendix.
\begin{lem}\label{lem:sup_tails} For any $\overline{\mathcal{C}}>0$, $t\in[0,\infty)$ and $v\in\R$,
\begin{align}
   \lim_{\mathcal{T}\to\infty}\sup_{u>0} \pk{\sup_{-\infty < s < -\mathcal{T}}
   \biggl(\frac{\overline{\mathcal{C}}}{\sigma(\delta(u))}\bigl(J(\delta(u)t) - J(\delta(u)s)\bigr) - (t-s)\biggr)>v}=0,\label{sup_tail}
\end{align}
and for any $j\in\sprod{n}$,
\begin{align}
   {\lim_{\mathcal{T}\to\infty}\pk{\sup_{-\infty < s < -\mathcal{T}}
   \bigl(\overline{\mathcal{C}}\bigl( B_{\lambda,j}(t) -B_{\lambda,j}(s)\bigr) - (t - s\bigr)\bigr)>v}=0.}\label{sup_tail_fBm}
\end{align}
\end{lem}

We fix some $\varepsilon>0$. Due to {Lemma \ref{lem:sup_tails}}, we can select a value of $T_{\varepsilon}>0$
such that simultaneously, for any $i\in\sprod{m}$, $j\in\sprod{n}$ {and $t_i<T_{\varepsilon}$},
\begin{align*}
    &\pk{\sup_{-\infty < s < -T_{\varepsilon}}\bigl({\mathcal{C}_j}\bigl( B_{\lambda,j}(t_i) -B_{\lambda,j}(s)\bigr) - (t_i - s\bigr)\bigr)>v_{ij}}<\varepsilon,\\
    &\pk{\sup_{-\infty < s < -T_{\varepsilon}}\biggl(\frac{{\mathcal{C}}_j}{\sigma(\delta(r_j(u)))}\bigl(J(\delta(r_j(u))t_i) - J(\delta(r_j(u))s)\bigr) - (t_i-s)\biggr)>v_{ij}}<\varepsilon.
\end{align*}
Hence, recalling the definition of $J_u^{\delta}(\cdot)$ from \eqref{J_delta_def},
\begin{align*}
    &\hspace{-1cm}\abs{ \,\pk{\begin{matrix}\overline{\vk X}_u^{\delta}(t_1)  >  \vk v_1 \\  \ldots \\ \overline{\vk X}_u^{\delta}(t_m)  >  \vk v_m \end{matrix}}- \pk{ \begin{matrix}{\vk{\mathfrak{X}}_{\lambda}(t_1)} > \vk v_1 \\ \ldots \\ {\vk{\mathfrak{X}}_{\lambda}(t_m)} > \vk v_m \end{matrix}}} \\
   \leqslant & \sum_{i=1}^{m}\sum_{j=1}^{n}\pk{\sup_{-\infty < s < -T_{\varepsilon}}\bigl({\mathcal{C}_j}\bigl( B_{\lambda,j}(t_i) -{B}_{\lambda,j}(s)\bigr) - (t_i - s\bigr)\bigr)>v_{ij}}\:+\\
    &\sum_{i=1}^{m}\sum_{j=1}^{n}\pk{\sup_{-\infty < s < -T_{\varepsilon}}\biggl(\frac{{\mathcal{C}}_j}{\sigma(\delta(r_j(u)))}\bigl(J(\delta(r_j(u))t_i) - J(\delta(r_j(u))s)\bigr) - (t_i-s)\biggr)>v_{ij}}\:+\\
    &\left|\:\pk{\begin{matrix}\sup\limits_{ -T_{\varepsilon}\leqslant s\leqslant 0}\bigl(\vk{\mathcal{C}}\circ(\vk J_u^{\delta}(t_1)-\vk J_u^{\delta}(s)) - (t_1-s)\vk 1\bigr)>\vk v_1 \\  \ldots \\ \sup\limits_{-T_{\varepsilon}\leqslant s\leqslant 0}\bigl(\vk{\mathcal{C}}\circ(\vk J^{\delta}_u(t_m)-\vk J_u^{\delta}(s)) - (t_m-s)\vk 1\bigr)>\vk v_m \end{matrix}}\right. \:-\\
    & \qquad\qquad\qquad\left.\pk{ \begin{matrix}\sup\limits_{-T_{\varepsilon} \leqslant s \leqslant 0}\bigl(\vk{\mathcal{C}}\circ(\vk B_{\lambda}(t_1)-\vk B_{\lambda}(s)) - (t_1-s)\vk 1\bigr) >\vk v_1 \\ \ldots \\ \sup\limits_{-T_{\varepsilon} \leqslant s \leqslant 0}\bigl(\vk{\mathcal{C}}\circ(\vk B_{\lambda}(t_m)-\vk B_{\lambda}(s)) - (t_m-s)\vk 1\bigr)>\vk v_m \end{matrix}}\:\right|.
\end{align*}
Using the continuous mapping theorem, according to \nelem{gauss_convergence} the last term converges to zero as $u\to\infty$. Hence,
\begin{align*}
    \limsup_{u\to 0}\abs{\pk{\begin{matrix}\overline{\vk X}_u^{\delta}(t_1)  >  \vk v_1 \\  \ldots \\ \overline{\vk X}_u^{\delta}(t_m)  >  \vk v_m \end{matrix}} - \pk{ \begin{matrix}{\vk{\mathfrak{X}}_{\lambda}(t_1)} > \vk v_1 \\ \ldots \\ {\vk{\mathfrak{X}}_{\lambda}(t_m)} > \vk v_m \end{matrix}}} \leqslant 2mn\varepsilon.
\end{align*}
Finally,   by sending $\varepsilon\downarrow 0$, we obtain the claimed assertion.
\QED

\subsection{Tightness}
{To verify that the sequence $\overline{\vk X}_u^{\delta}(\cdot)$  is tight in $C([0,\infty)\to\R^n)$,
it is sufficient to prove that sequence $\overline{\vk X}_u^{\delta}(t)$ for $t\in[0,T]$ is tight in $C([0,T]\to\R^n)$
for every $T>0$. Hence, due to \cite[Theorem 7.5]{billingsley2013convergence} and taking into account \eqref{sup_fdd_claim},
it is enough}
to show that, for any $i\in\sprod{n}$, any $T>0$ {and any $\eta>0$},
\begin{align}
    \lim_{\zeta \downarrow 0}\limsup_{u\to\infty}\pk{\sup_{(s,t)\in{\mathscr S}_{T,\zeta}}\abs{\frac{\overline{X}_{u,i}\bigl(\delta(r_i(u))\,t\bigr)}{\sigma\bigl(\delta(r_i(u))\bigr)} - \frac{\overline{X}_{u,i}\bigl(\delta(r_i(u))\,s\bigr)}{\sigma\bigl(\delta(r_i(u))\bigr)}}\geqslant\eta} = 0,\label{uniformly_equicontinuous_component_bar}
\end{align}
where ${\mathscr S}_{T,\zeta}:=\{s,t\in[0,T]: |t-s|\leqslant \zeta\}.$
Using \eqref{X_bar_def}, {we obtain that
\begin{align*}
   { \overline{X}_{u,i}\bigl(\delta(r_i(u))\,t\bigr)}&=\sup_{-\infty < s < t}\biggl(\mathcal{C}_i\bigl(J\bigl(\delta(r_i(u))\,t\bigr)-J\bigl(\delta(r_i(u))\,s\bigr)\bigr) - r_i(u)\,\delta(r_i(u))\,(t-s)\biggr)\\
    &=\,\mathcal{C}_i\sup_{-\infty < s < t}\biggl(J\bigl(\delta(r_i(u))\,t\bigr)-J\bigl(\delta(r_i(u))\,s\bigr)\bigr) - r_i(u)\,\delta(r_i(u))\,\frac{1}{\mathcal{C}_i}(t-s)\biggr)\\
    &= \,\mathcal{C}_i\,\overline{X}^{\star}_{r_i(u),1/\mathcal{C_i}}\bigl(\delta(r_i(u))\,t\bigr),
\end{align*}
where}
\begin{align*}
    \overline{X}_{{u,c}}^{\star}(t) := \sup_{-\infty < s < t}\bigl(\bigl(J(t) -J(s)\bigr) - {u\,c\,(t-s)}\bigr) \geqslant 0.
\end{align*}
{Hence, to verify \eqref{uniformly_equicontinuous_component_bar} it is sufficient to show that, for any $c>0$, $\eta>0$,}
\begin{align}
    \lim_{\zeta \downarrow 0}\limsup_{u\to\infty}\pk{\sup_{(s,t)\in{\mathscr S}_{T,\zeta}}\abs{\frac{\overline{X}_{{u,c}}^{\star}\bigl(\delta(u)t)\bigr)}{\sigma\bigl(\delta(u)\bigr)} - \frac{\overline{X}_{{u,c}}^{\star}\bigl(\delta(u)s\bigr)}{\sigma\bigl(\delta(u)\bigr)}}\geqslant\eta} = 0.\label{uniformly_equicontinuous_bar_2}
\end{align}

Recall a standard representation of the reflected process:  we can write, for any $t>s$,
\begin{align*}
    \overline{X}_{{u,c}}^{\star}(t) = \overline{X}_{{u,c}}^{\star}(s) + J(t)-J(s) - {u\,c\,(t-s)} +
    \max\biggl(0, \sup_{s<v<t}\bigl(-\overline{X}_{{u,c}}^{\star}(s) - (J(v)-J(s))+{u\,c\,(v-s)}\bigr)\biggr);
\end{align*}
see for instance \cite[\S 2.4]{DeM15}. Hence, for all $u$, and for any given $T,\zeta>0$,
\begin{align}\notag
\sup_{(s,t)\in{\mathscr S}_{T,\zeta}}    \abs{\overline{X}_{{u,c}}^{\star}(\delta(u)\,t) - \overline{X}_{{u,c}}^{\star}(\delta(u)\,s)}
&=
  \sup_{(s,t)\in{\mathscr S}_{T,\zeta}}  \biggl|J\bigl(\delta(u)\,t\bigr) -J\bigl(\delta(u)\,s\bigr)- { u\,c\,\delta(u)}\,(t-s)\\
    &\:\:\:\:\hspace{-3.5cm}+ \max\biggl(0, \sup_{{s<v<t}}\bigl(-\overline{X}_{{u,c}}^{\star}(\delta(u)\,s) - (J({\delta(u)\,v})-J({\delta(u)\,s}))+{u\,c\,\delta(u)\,(v-s)}\bigr)\biggr)\biggr|\notag\nonumber\\
   &\:\:\:\: \hspace{-3.9cm}\leqslant
  \sup_{(s,t)\in{\mathscr S}_{T,\zeta}}  \biggl|J\bigl(\delta(u)\,t\bigr) -J\bigl(\delta(u)\,s\bigr)- {u\,c\,\delta(u)}\,(t-s)\biggr|\nonumber\\
    &\:\:\:\: \hspace{-3.5cm} +
    \sup_{(s,t)\in{\mathscr S}_{T,\zeta}}
    \max\biggl(0, \sup_{{s<v<t}}\bigl(-\overline{X}_{{u,c}}^{\star}(\delta(u)\,s) - (J({\delta(u)\,v})-J({\delta(u)\,s}))+{u\,c\,\delta(u)}\,(v-s)\bigr)\biggr)\notag\nonumber\\
     &{\:\:\:\:\hspace{-3.9cm}\leqslant
  \sup_{(s,t)\in{\mathscr S}_{T,\zeta}}  \biggl|J\bigl(\delta(u)\,t\bigr) -J\bigl(\delta(u)\,s\bigr)- {u\,c\,\delta(u)}\,(t-s)\biggr|\nonumber}\\
    &\:\:\:\:\:\:\:\: \hspace{-3.5cm}{+
    \max\biggl(0, \sup_{{\substack{(s,t)\in{\mathscr S}_{T,\zeta}\\ s<v<t}}}\bigl( - (J({\delta(u)\,v})-J({\delta(u)\,s}))+{u\,c\,\delta(u)}\,(v-s)\bigr)\biggr)\notag\nonumber}\\
     &\:\:\:\: \hspace{-3.9cm}=
      \sup_{(s,t)\in{\mathscr S}_{T,\zeta}}
      2\abs{J\bigl(\delta(u)\,t\bigr) -J\bigl(\delta(u)\,s\bigr)- {u\,c\,\delta(u)}\,(t-s)},
    \label{X_star_bound_in_J}
\end{align}
{where (a)~in the first inequality we use the triangle inequality, (b)~in the second inequality the fact that $\overline{X}_{{u,c}}^{\star}(\delta(u)\,s)\geqslant 0$, and (c)~in the last equality the fact that $-J(\cdot)$ has the same distribution as $J(\cdot)$.} Applying the upper bound \eqref{X_star_bound_in_J} in \eqref{uniformly_equicontinuous_bar_2}, {and recalling that by \eqref{delta_def}} we have that $u\,\delta(u) = \sigma(\delta(u))$, we obtain after elementary manipulations, for any value of $u$,
\begin{align*}
    \pk{\sup_{(s,t)\in{\mathscr S}_{T,\zeta}}
    \abs{\frac{\overline{X}_{{u,c}}^{\star}\bigl(\delta(u)t)\bigr)}{\sigma\bigl(\delta(u)\bigr)} -
    \frac{\overline{X}_{{u,c}}^{\star}\bigl(\delta(u)s\bigr)}{\sigma\bigl(\delta(u)\bigr)}}\geqslant\eta}
    \leqslant
    \pk{\sup_{(s,t)\in{\mathscr S}_{T,\zeta}}
       \abs{\frac{J\bigl(\delta(u)t)\bigr)-J\bigl(\delta(u)s)\bigr)}{\sigma\bigl(\delta(u)\bigr)} - {c}\,(t-s)}\geqslant\frac{\eta}{2}}.
\end{align*}
We now consider the right hand side of the previous display in the regime that $u$ grows large. Due to \nelem{gauss_convergence}, {in combination with the continuous mapping theorem,} we obtain that, as $u\to\infty$,
\begin{align*}
    \pk{\sup_{(s,t)\in{\mathscr S}_{T,\zeta}}
       \abs{\frac{J\bigl(\delta(u)t)\bigr)-J\bigl(\delta(u)s)\bigr)}{\sigma\bigl(\delta(u)\bigr)} - {c}\,(t-s)}\geqslant\frac{\eta}{2}}
    \to
    \pk{\sup_{(s,t)\in{\mathscr S}_{T,\zeta}}\abs{B_{\lambda}(s,t) - {c}\,(t-s)}\geqslant\frac{\eta}{2}};
\end{align*}
here $B_{\lambda}(s,t):=B_\lambda(t)-B_\lambda(s)$ with $B_\lambda(\cdot)$ denoting a fractional Brownian motion with Hurst parameter $\lambda$. Hence, for sufficiently large enough  values of $u$,
\begin{align*}
    \pk{\sup_{(s,t)\in{\mathscr S}_{T,\zeta}}\abs{\frac{\overline{X}_{{u,c}}^{\star}\bigl(\delta(u)t)\bigr)}
          {\sigma\bigl(\delta(u)\bigr)} - \frac{\overline{X}_{{u,c}}^{\star}\bigl(\delta(u)s\bigr)}{\sigma\bigl(\delta(u)\bigr)}}
          \geqslant\eta}
    \leqslant
    2\,
\pk{\sup_{(s,t)\in{\mathscr S}_{T,\zeta}}\abs{B_{\lambda}(s,t) - {c}\,(t-s)}\geqslant\frac{\eta}{2}}.
\end{align*}
Now \eqref{uniformly_equicontinuous_bar_2} follows from the fact that, {for any $c>0$ and $\eta>0$}
\begin{align}
    \lim_{\zeta \downarrow 0}
    \pk{\sup_{(s,t)\in{\mathscr S}_{T,\zeta}}\abs{B_{\lambda}(s,t) - {c}\,(t-s)}\geqslant\frac{\eta}{2}}=0,\label{uniformly_equicontinuous_fBm}
\end{align}
which completes the proof.
We have included the proof of \eqref{uniformly_equicontinuous_fBm} in the appendix.
\QED

\section{Discussion and concluding remarks} We have analyzed a Gaussian-driven queueing network with a routing mechanism that prohibits both the merging of multiple streams and the formation of loops. The primary results focus on functional limit theorems for the joint workload process, particularly in the light-traffic and heavy-traffic regimes.

Our work has intriguing connections to the recent paper \cite{DAI}, which examines a generalized Jackson network with a general routing matrix. Based on the results presented in our paper, one anticipates specific possible extensions of the results in \cite{DAI}. First, while \cite{DAI} focuses on limit theorems for the stationary distribution of the queue lengths, it raises the question of whether these results extend to functional versions. In the limiting regime considered, one would anticipate that the joint workload process converges to independent reflected Brownian motions. Second, it would be valuable to explore cases where multiple queues are essentially equally heavily loaded, as we would expect such groups of workloads to converge to {\it dependent} reflected Brownian motions.

Another avenue for potential extensions focuses on the Gaussian framework presented in this paper, but with a more general routing mechanism in place. The approach outlined here relies heavily on the representation in \eqref{Q_vk_def}, which is specific to the current routing structure and does not extend to a general routing matrix. Thus, a completely new approach would be required to address this general class of models.
Alternatively, one could retain the routing matrices considered in the present paper but work with Lévy inputs. It is anticipated that this analysis could be feasible, given the availability of explicit results for the transform of the joint workload \cite{DDR,DeM15}.

\appendix

\section{Proofs of lemmas}\label{Section:Appendix}

\prooflem{lem:sup_tails} {We consider only \eqref{sup_tail}; \eqref{sup_tail_fBm} can be obtained by a similar reasoning.} Let $\overline{\mathcal{C}}>0$, $t\in[0,\infty)$ and $v\in\R$ be given. Using that the function under the limit in \eqref{sup_tail} is non-increasing in $\mathcal{T}$, we can restrict ourselves to $\mathcal{T}\in\N$ and $\mathcal{T}>2v-t$, and $u$ large enough such that $\delta(u)<\mathfrak{a}$, where $\mathfrak{a}$ is as defined in \nelem{lem:F_uniform_bound}. For any $i\in\N$, $s_1,s_2\in[-i-1,-i]$, by applying \eqref{F_uniform_bound}, we find that
\begin{align*}
    &\E{\left(\frac{2\overline{\mathcal{C}}J\bigl(\delta(u)t\bigr) - 2\overline{\mathcal{C}}J\bigl(\delta(u)s_1\bigr)}{(t-s_1)\sigma\bigl(\delta(u)\bigr)}-\frac{2\overline{\mathcal{C}}J\bigl(\delta(u)t\bigr) - 2\overline{\mathcal{C}}J\bigl(\delta(u)s_2\bigr)}{(t-s_2)\sigma\bigl(\delta(u)\bigr)}\right)^2}\\
    &\qquad\leqslant \frac{4\overline{\mathcal{C}}^2(s_1-s_2)^2\sigma^2\bigl(\delta(u)t\bigr)}{(t-s_1)^2(t-s_2)^2\sigma^2\bigl(\delta(u)\bigr)}+\frac{4\overline{\mathcal{C}}^2(t-s_1)^2\sigma^2\bigl(\delta(u)\abs{s_1-s_2}\bigr)}{(t-s_1)^2(t-s_2)^2\sigma^2\bigl(\delta(u)\bigr)}\\
    &\qquad\qquad+\frac{4\overline{\mathcal{C}}^2(s_1-s_2)^2\sigma^2\bigl(\delta(u)s_1\bigr)}{(t-s_1)^2(t-s_2)^2\sigma^2\bigl(\delta(u)\bigr)}\\
    &\qquad\leqslant \frac{4\overline{\mathcal{C}}^2\mathfrak{C}^2(t^{\gamma_{\infty}}+t^{\gamma_{0}})^2(s_1-s_2)^2}{(t-s_1)^2(t-s_2)^2}+\frac{4\overline{\mathcal{C}}^2\mathfrak{C}^2(t-s_1)^2\left(\abs{s_1-s_2}^{\gamma_{\infty}}+\abs{s_1-s_2}^{\gamma_{0}}\right)^2}{(t-s_1)^2(t-s_2)^2}\\
    &\qquad\qquad+\frac{4\overline{\mathcal{C}}^2\mathfrak{C}^2(\abs{s_1}^{\gamma_{\infty}}+\abs{s_1}^{\gamma_{0}})^2(s_1-s_2)^2}{(t-s_1)^2(t-s_2)^2}\\
    &\qquad\leqslant 4\overline{\mathcal{C}}^2\mathfrak{C}^2(s_1-s_2)^2+8\overline{\mathcal{C}}^2\mathfrak{C}^2\abs{s_1-s_2}^{2\min(\gamma_{\infty},\gamma_{0})}+4\overline{\mathcal{C}}^2\mathfrak{C}^2(s_1-s_2)^2\\
    &\qquad\leqslant 16\overline{\mathcal{C}}^2\mathfrak{C}^2\abs{s_1-s_2}^{2\min(\gamma_{\infty},\gamma_{0})}.
\end{align*}

Applying Piterbarg's inequality \cite[Thm.\ 8.1]{piterbarg1996asymptotic}, we obtain that, for some constant $C>0$,
\begin{align}
    \varphi({\mathcal T},u):=&\,\pk{\sup_{-\infty < s < -\mathcal{T}}\frac{\overline{\mathcal{C}}}{\sigma\bigl(\delta(u)\bigr)}\biggl(\bigl({J\bigl(\delta(u)t\bigr) - J\bigl(\delta(u)s\bigr)}\bigr) - (t-s)\biggr)>v}\notag\\
    \leqslant&\, \sum_{i=\mathcal{T}}^{\infty}\pk{\sup_{s\in[-i-1,-i]}\frac{2\overline{\mathcal{C}}}{\sigma\bigl(\delta(u)\bigr)}\left(\frac{J\bigl(\delta(u)t\bigr) - J\bigl(\delta(u)s\bigr)}{(t-s)}\right)>1}\leqslant\sum_{i=\mathcal{T}}^{\infty}\frac{C\sigma^{\star\star}_{u}(i)}{\sqrt{2\pi}}e^{-\frac{1}{2\left(\sigma^{\star\star}_{u}(i)\right)^2}},\label{Pit_for_sup}
\end{align}
where
\begin{align*}
    \bigl(\sigma^{\star\star}_{u}(i)\bigr)^2 &= \sup_{s\in[-i-1,-i]}\E{\left(\frac{2\overline{\mathcal{C}}}{\sigma\bigl(\delta(u)\bigr)}\right)^2\left(\frac{J\bigl(\delta(u)t\bigr) - J\bigl(\delta(u)s\bigr)}{(t-s)}\right)^2}\\ &= 4\overline{\mathcal{C}}^2\sup_{s\in[-i-1,-i]}\frac{\sigma^2\bigl(\delta(u)(t-s)\bigr)}{(t-s)^2\sigma^2\bigl(\delta(u)\bigr)}= 4\overline{\mathcal{C}}^2\sup_{s\in[-i-1,-i]}\frac{\mathcal{F}^2_{\delta(u)}(t-s)}{(t-s)^2}.
\end{align*}
According to \eqref{F_uniform_bound_t>1} (bearing in mind that $t-s\geqslant t+\mathcal{T}>1$), for $u$ sufficiently large,
\begin{align}
    \sigma^{\star\star}_{u}(i) \leqslant \frac{4\,\mathfrak{C}\,\overline{\mathcal{C}}}{(i+t)^{1-\gamma_{\infty}}}.\label{sigma_sup_bound}
\end{align}
Hence, combining \eqref{Pit_for_sup} with \eqref{sigma_sup_bound}, for $\mathcal{T}$ sufficiently large, we have the upper bound
\begin{align*}
    \varphi({\mathcal T},u) &\leqslant \sum_{i=\mathcal{T}}^{\infty}\psi(i)\leqslant \int_{\mathcal{T}-1}^{\infty} \psi(x)\,\td x,\quad\:\:\psi(x):=\frac{4\,C\,\mathfrak{C}\,\overline{\mathcal{C}}}{\sqrt{2\pi}(x+t)^{1-\gamma_{\infty}}}\exp\left({-\frac{(x+t)^{2-2\gamma_{\infty}}}{32\mathfrak{C}^2\overline{\mathcal{C}}^2}}\right);
\end{align*}
importantly, this upper bound is uniform in $u$.
Now \eqref{sup_tail} directly follows from   $\int_{0}^{\infty} \psi(x)\,\td x<\infty$,
where it has been used that $\gamma_{\infty}<1$.
\QED

\medskip

\prooflem{lem:F_uniform_bound} Fix $\varepsilon>0$ such that $\varepsilon<\min(\lambda, 1-\lambda)$. Due to \textbf{L1}, there exist $\mathfrak{a}>0$ and constants $C_1,C_2>0$ such that for any $x,t$ if $x\leqslant \mathfrak{a}$ and $tx\leqslant \mathfrak{a}$, we have
\begin{align}
    \mathcal{F}_{x}(t)\leqslant C_1 t^{\lambda-\varepsilon},\label{F_uniform_bound_part_1}
\end{align}
while for any $x\leqslant\mathfrak{a}$
\begin{align}
    \sigma(x)>C_2x^{\lambda+\varepsilon}.\label{sigma>}
\end{align}
Using that the function $\sigma(x)$ is continuous, by \textbf{L2} we can find a constant $C_3$ such that for all $x>\mathfrak{a}$
\begin{align}
    \sigma(x)\leqslant C_3 x^{\beta}.\label{sigma<}
\end{align}
Hence, for $x,t$ such that $x\leqslant\mathfrak{a}$ and $tx>\mathfrak{a}$, upon combining \eqref{sigma>} and \eqref{sigma<},
\begin{align}
    \mathcal{F}_{x}(t) = \frac{\sigma(tx)}{\sigma(x)}\leqslant \frac{C_3(tx)^{\beta}}{C_2 x^{\lambda+\varepsilon}} = C_4t^{\beta}x^{\beta-\lambda-\varepsilon} .\label{sigma_frac_bound*}
\end{align}
In case $\beta-\lambda-\varepsilon\geqslant 0$ we have from \eqref{sigma_frac_bound*}, due to $x\leqslant\mathfrak{a}$,
\begin{align}
    \mathcal{F}_{x}(t)\leqslant C_4\mathfrak{a}^{\beta-\lambda-\varepsilon}t^{\beta}\label{tx>a_bound_1}.
\end{align}
In case $\beta-\lambda-\varepsilon< 0$ we have from \eqref{sigma_frac_bound*}, due to $tx>\mathfrak{a}$,
\begin{align}
    \mathcal{F}_{x}(t)\leqslant C_4(tx)^{\beta-\lambda-\varepsilon} t^{\lambda+\varepsilon} \leqslant C_4\mathfrak{a}^{\beta-\lambda-\varepsilon} t^{\lambda+\varepsilon}.\label{tx>a_bound_2}
\end{align}
Combining \eqref{tx>a_bound_1} with \eqref{tx>a_bound_2}, we obtain that there exists a constant $C_5>0$ such that if $x\leqslant\mathfrak{a}$ and $tx>\mathfrak{a}$, then
\begin{align}
   \mathcal{F}_{x}(t)\leqslant C_5 t^{\max(\beta,\lambda+\varepsilon)}.\label{F_uniform_bound_part_2}
\end{align}
Hence, combining \eqref{F_uniform_bound_part_1} with \eqref{F_uniform_bound_part_2}, both claims \eqref{F_uniform_bound} and \eqref{F_uniform_bound_t>1} follow by picking $\mathfrak{C} :=\max(C_1,C_5)$, $\gamma_{\infty} := \max\bigl(\max(\beta,\lambda+\varepsilon), \lambda-\varepsilon\bigr)$, and $\gamma_{0} := \min\bigl(\max(\beta,\lambda+\varepsilon), \lambda-\varepsilon\bigr)$.
\QED

\medskip

\proofrem{rem:F_uniform_bound_heavy} Fix $\varepsilon>0$ such that $\alpha-\varepsilon>0$ and $\alpha+\varepsilon<1$. Due to \textbf{H1}, there exist $\mathfrak{a}>0$ and constant $C_1>0$ such that for all $x,t$ if $x\geqslant \mathfrak{a}$ and $tx\geqslant \mathfrak{a}$, then
\begin{align}
    \mathcal{F}_{x}(t)\leqslant C_1 t^{\lambda+\varepsilon}.\label{F_uniform_bound_heavy_infty}
\end{align}
Following the same arguments as in the proof of \nelem{lem:F_uniform_bound}, we obtain that there exist constant $C_2>0$ such that for all $x,t>0$ if $x\geqslant \mathfrak{a}$ and $tx\leqslant \mathfrak{a}$, then
\begin{align}
   \mathcal{F}_{x}(t)\leqslant C_2 t^{\min(\beta,\alpha-\varepsilon)}.\label{F_uniform_bound_heavy_0}
\end{align}
Then, combining \eqref{F_uniform_bound_heavy_infty} and \eqref{F_uniform_bound_heavy_0} we obtain that all $x,t>0$ if $x\geqslant \mathfrak{a}$ then
\begin{align*}
    \mathcal{F}_{x}(t)\leqslant (C_1+C_2)\bigl(t^{\lambda+\varepsilon} + t^{\min(\beta,\alpha-\varepsilon)}\bigr),
\end{align*}
and if $t>1$ then due to \eqref{F_uniform_bound_heavy_infty}
\begin{align*}
    \mathcal{F}_{x}(t)\leqslant C_1 t^{\lambda+\varepsilon}.
\end{align*}
Hence, the claim follows.
\QED

\medskip

{\textbf{Proof of Equation} \eqref{uniformly_equicontinuous_fBm}: We prove this equation by appealing to Piterbarg's inequality, as follows. For any $\lambda\in(0,1)$, $c>0$, $\zeta>0$ and $\eta>0$, by the triangle inequality,
\begin{align*}
    \pk{\sup_{(s,t)\in{\mathscr S}_{T,\zeta}}\abs{B_{\lambda}(s,t) - {c}\,(t-s)}\geqslant\frac{\eta}{2}}&\leqslant\pk{\sup_{(s,t)\in{\mathscr S}_{T,\zeta}}\abs{B_{\lambda}(s,t)} + \sup_{(s,t)\in{\mathscr S}_{T,\zeta}}\abs{c(t-s)}\geqslant\frac{\eta}{2}}\\
    &\leqslant\pk{\sup_{(s,t)\in{\mathscr S}_{T,\zeta}}\abs{B_{\lambda}(s,t)}\geqslant\frac{\eta}{2} - c\zeta},
\end{align*}
where we recall that $B_\lambda(s,t):=B_\lambda(t)-B_\lambda(s)$ and  ${\mathscr S}_{T,\zeta}:=\{s,t\in[0,T]: |t-s|\leqslant \zeta\}.$
Hence, for any $\zeta < \eta/(4c)$,
\begin{align*}
    \pk{\sup_{(s,t)\in{\mathscr S}_{T,\zeta}}\abs{B_{\lambda}(s,t) - {c}\,(t-s)}\geqslant\frac{\eta}{2}}&\leqslant\pk{\sup_{(s,t)\in{\mathscr S}_{T,\zeta}}\abs{B_{\lambda}(s,t)}\geqslant\frac{\eta}{4}}\leqslant 2\,\pk{\sup_{(s,t)\in{\mathscr S}_{T,\zeta}}B_{\lambda}(s,t)\geqslant\frac{\eta}{4}},
\end{align*}
where in the second inequality we use the union bound together with the fact that $B_{\lambda}(\cdot)$ and $-B_{\lambda}(\cdot)$ have the same distribution. Noticing that, for any ${s_1,s_2,t_1,t_2\in[0,T]}$,
\begin{align*}
    \E{\bigl(B_{\lambda}(s_1,t_1)\bigr)^2} &= \abs{t_1-s_1}^{2\lambda},\\
    \E{\bigl(B_{\lambda}(s_1,t_1) - B_{\lambda}(s_2,t_2)\bigr)^2} &\leqslant 2\left(\E{\bigl(B_{\lambda}(s_1,s_2)\bigr)^2} + \E{\bigl(B_{\lambda}(t_1,t_2)\bigr)^2}\right)\\
    &\leqslant 2^{2-\lambda}\left(\abs{s_1-s_2}^2 + \abs{t_1-t_2}^2\right)^{\lambda}
\end{align*}(where in the last line we apply the generalized mean inequality),
due to Piterbarg's inequality \cite[Theorem 8.1]{piterbarg1996asymptotic},
there exists a constant $C>0$ such that
\begin{align*}
    \pk{\sup_{(s,t)\in{\mathscr S}_{T,\zeta}} B_{\lambda}(s,t) \geqslant\frac{\eta}{4}}&\leqslant {CT^{2}\left(\frac{\eta}{4}\right)^{4/(2\lambda)}\frac{\sigma_{\star}}{\sqrt{2\pi}(\eta/4)}}\exp\left({-\frac{(\eta/4)^2}{2\sigma_{\star}^{2}}}\right),\\
    \sigma_{\star}^{2} &:= \sup_{(s,t)\in{\mathscr S}_{T,\zeta}}\E{(B_{\lambda}(s,t) )^2} \leqslant \zeta^{2\lambda}.
\end{align*}
Hence, for any $\lambda\in(0,1)$, $c>0$, and $\eta>0$,
\begin{align*}
    \lim_{\zeta\downarrow 0}\pk{\sup_{(s,t)\in{\mathscr S}_{T,\zeta}}\abs{B_{\lambda}(s,t) - {c}\,(t-s)}\geqslant\frac{\eta}{2}}&\leqslant\frac{{4}\,CT^2\eta^{{2/\lambda}}}{2^{{4/\lambda}}\sqrt{2\pi}\,\eta}\lim_{\zeta\downarrow 0}\zeta^{\lambda}\exp\left({-\frac{(\eta/4)^2}{2\zeta^{2\lambda}}} \right)= 0,
\end{align*}
which completes the proof.}
\QED

{\small
\bibliographystyle{ieeetr}
\bibliography{EEEA}}
\end{document}